\documentclass{amsart}

\usepackage{amssymb}
\usepackage{amsmath}
\usepackage{amsthm}
\usepackage{amscd}

\newcommand{\Z}{{\mathbb Z}}            % the integers Z
\newcommand{\Q}{{\mathbb Q}}            % the rationals Q
\newcommand{\R}{{\mathbb R}}            % the reals R
\newcommand{\C}{{\mathbb C}}            % the complexes C
\newcommand{\Ocal}{{\mathcal O}}        % the ring of integers
\renewcommand{\c}{{\mathbf c}}            % normalized capacity
\renewcommand{\d}{{\mathbf d}}            % normalized nth diameter
\newcommand{\g}{{\mathbf g}}            % Adelic generalized Green's function
\newcommand{\Fp}{{\mathbf F_p}}         % finite field of p elements
\newcommand{\Kbar}{{\bar K}}            % alg. closure of K
\newcommand{\Lbar}{{\bar L}}            % alg. closure of L
\newcommand{\Kvbar}{{\bar K_v}}         % alg. closure of K_v 
\newcommand{\Qbar}{{\overline{\Q}}}
\newcommand{\M}{{\mathfrak{M}}}         % mandelbrot set

\newcommand{\into}{\hookrightarrow}     % an injection into
\renewcommand{\P}{{\mathbf P}}          % Projective space

\newcommand{\F}{{\mathbf F}}                    % Used for finite fields
\newcommand{\FF}{{\mathbb F}}                    % Used for finite fields

\newcommand{\hhat}{\hat{h}}             % canonical height
\newtheorem{theorem}{Theorem}[section]

\newtheorem{cor}[theorem]{Corollary}
\newtheorem{lemma}[theorem]{Lemma}
\newtheorem{prop}[theorem]{Proposition}

\theoremstyle{remark}
\newtheorem{remark}[theorem]{{\bf Remark}}

\newtheorem{example}[theorem]{{\rm\bf  Example}}

\newcounter{listcounter1}
\newenvironment{list1}{
  \begin{list}{(\arabic{listcounter1}).\hfill}{
    \usecounter{listcounter1}
    \setlength{\leftmargin}{18pt}
    \setlength{\labelwidth}{18pt}
    \setlength{\labelsep}{0pt}
    \setlength{\topsep}{0pt}
  }
}{
  \end{list}
}


\begin{document}

\title[Canonical heights and polynomial dynamics]
{Canonical heights, transfinite diameters, and polynomial dynamics}

% \date{First draft September 2002; Latest revision March 2003}

\author[Matthew Baker]{Matthew H. Baker}
\author[Liang-Chung Hsia]{Liang-Chung Hsia}
\email{mbaker@math.uga.edu \\
       hsia@math.ncu.edu.tw}
\address{Department of Mathematics,
         University of Georgia, Athens, GA 30602-7403, USA}
\address{Department of Mathematics, 
         National Central University, 
         Chung-Li, 32054, Taiwan}

\thanks{The first author's research was supported in part
by an NSF Postdoctoral Research Fellowship.}
\thanks{The second author's research was supported by an NSC research
  grant 91-2115-M-008-010.}
\thanks{The authors would like to thank Robert Rumely and Joseph
  Silverman for helpful conversations. We would also like to thank the 
  referee for pointing out some errors in the original version of this
  paper.}

\begin{abstract}
Let $\phi(z)$ be a polynomial of degree at least 2
with coefficients in a number field $K$.
Iterating $\phi$ gives rise to a dynamical system and a corresponding
canonical height function $\hhat_\phi$, as defined by Call and
Silverman.  We prove a simple product formula
relating the transfinite diameters of
the filled Julia sets of $\phi$ over various completions of $K$,
and we apply this formula to give a generalization of
Bilu's equidistribution theorem
for sequences of points whose canonical heights tend to zero.
\end{abstract}

\maketitle
%-----------------------------------------------------------

\section{Introduction}

There are a number of interesting results in arithmetic geometry and
in the theory of complex dynamical systems which deal with
equidistribution properties of sequences of ``small'' points.  Examples on the arithmetic side include
Bilu's equidistribution theorem (see \cite{Bilu})
and the equidistribution theorem of Szpiro, Ullmo, and Zhang (see
\cite{SUZ}), which imply the ``generalized Bogomolov conjecture'' for
subvarieties of algebraic tori and abelian varieties, respectively.
(See \cite{Bilu}, \cite{Ullmo}, and \cite{Zhang} for further details).
On the complex dynamics side, there is for example an important
equidistribution result for the backward iterates of a point under a polynomial
due to Brolin (see \cite{Brolin} or  \cite[Theorem 6.5.8 ]{Ransford}).
Brolin's result is closely related to both the strengths and
weaknesses of the so-called ``inverse iteration method'' for drawing Julia sets
of polynomials on a computer.  (See \cite[Page 207]{Ransford} and
\cite[Chapter 2]{Peitgen-Richter} for a further discussion).

One of our goals in this paper is to simultaneously generalize Bilu's
theorem and Brolin's theorem, and to explore the interplay of both results
with potential theory.  The connection between Bilu's theorem and potential theory was already
pointed out by Rumely in \cite{Rumely} and by Bombieri in \cite{Bombieri}.
Because some of the sets
we deal with in this paper are rather exotic (they can be thought of
as ``$p$-adic fractals''), we have found it convenient to
base all of our proofs on the purely algebraic notion of transfinite
diameter, rather than on the equivalent potential-theoretic notion of capacity.  

Another goal of this paper is to promote the philosophy that
interesting information about dynamical systems can be found by
looking ``adelically'', rather than just over $\C$.  A simple but
instructive example is the observation that $0$ is not a periodic
point under iteration of the polynomial $z^2 - 3/2$.  This is not
obvious from plotting the iterates in the real line, but it is clear
if we look at what is happening $2$-adically.  More generally, it
seems useful to try to put the $p$-adic theory of polynomial iteration on
an equal footing with the complex theory \cite{Benedetto-Component, 
Benedetto-Hyperbolic, Hsia1, Hsia2, Juan}.  
%\helpme{Give some references to papers of Hsia, Benedetto, etc.}
However, complications inevitably arise from the fact that $\C_p$ is
not locally compact.  For example, in Section~\ref{thm1proofsection} we compute the transfinite diameter of
the filled Julia set of a polynomial $\phi$ over $\C_p$, but the computation is not
quite as simple as its complex analogue due to the failure of local compactness.  This calculation leads to a
simple global product formula, which in turn is closely related to the
equidistribution results alluded to above.

\medskip

We have aimed to make this paper as self-contained as possible.  
We hope that this makes the
paper a useful introduction to certain aspects of the theory of canonical heights, potential theory, 
and dynamical systems over local and global fields.

\section{Terminology and conventions}
\label{terminologysection}

Throughout this paper, $K$ will denote a {\it global field}, which for
our purposes will mean a field equipped with a product
formula and satisfying the Northcott finiteness property.\footnote{This
is not the standard definition of a global field.}  The main examples
are number fields and function fields of curves over finite fields.

To say that $K$ is equipped with a product formula means that we are given a family $M_K$
of absolute values $x \mapsto |x|_v$ on $K$ (where $v$ denotes the corresponding
additive valuation) and a family 
$\{ n_v \} \; (v \in M_K)$ of positive real numbers satisfying
the following properties (compare with \cite[Section 2.1]{SerreLMW}):
\begin{itemize}
\item All but finitely many $v \in M_K$ are ultrametric\footnote{We
    will use the terms ultrametric and non-archimedean interchangeably
    throughout this paper.}, meaning that 
$|x|_v = c^{-v(x)}$ for some real additive valuation (which we also
    denote by $v$) on $K$ and some $c>1$.
\item The absolute values which are not ultrametric are induced by embeddings
$K \into \C$.  
\item For all $x \in K^*$, we have $|x|_v = 1$ for all but finitely
  many $v \in M_K$.
\item If we set $\| \quad \|_v := |\quad |_v^{n_v}$, 
then for all $x \in K^*$, we have the product formula
\[
\prod_{v \in M_K} \| x \|_v = 1.
\]
\end{itemize}

Fix a place $v \in M_K$, and let $K_v$ denote the completion of $K$ with
respect to $|\quad|_v$.
If we fix an algebraic closure $\Kvbar$ of $K_v$, then it is well-known
(see  \cite[Section~II.1]{LangANT}) that $|\quad|_v$ 
extends in a unique way to an
absolute value on $\Kvbar$.  
Moreover, if we denote by $\C_v$ the completion of $\Kvbar$,
then $|\quad|_v$ extends uniquely to $\C_v$ by continuity.
For convenience, we often fix without comment a distinguished embedding of $\Kbar$ into $\C_v$ for each
$v \in M_K$.

\medskip

We say that $K$ satisfies the {\it Northcott finiteness property} if
for each integer $d \geq 1$ and each positive real number $M$, there
are only finitely many elements $\alpha \in
\Kbar$ of degree at most $d$ over $K$ and absolute logarithmic height at
most $M$.  
(See \cite[Section~2.2]{SerreLMW} or
Section~\ref{resultsection} 
below for a definition of the absolute logarithmic
height).  For example, it is well-known that 
number fields and function fields of  curves over finite fields 
(equipped with the standard product formula structure as
in \cite{SerreLMW}) have this property.
% \fixme{More example other than these two?
% How about let's concentrate on number fields and function 
% fields over finite field but mention that Theorem~\ref{thm1} is 
% valid for any global field, say?}

\medskip

Throughout this paper, $L$ will denote a valued field.
For ease of exposition, we will say that a valued field $L$ is a 
{\it local field} if $L$
is either $\R$ or $\C$ (with the usual absolute value) or a field which is 
complete with respect to an ultrametric absolute value. If 
the absolute value on $L$ arises from some $v\in M_K$ then the
absolute value will be denoted by $|\quad|_v$ as usual and the completion
of $\Lbar$ will be denoted by $\C_v$.
Otherwise we simply denote the absolute value on $L$ by $|\quad|$
(dropping the subscript $v$). For a non-archimedean field $L$,
the set $\{z\in L \; : \; |z| \le R\}$ will be called a closed disc (of radius
$R$) and the set $\{z\in L \; : \; |z| < R\}$ will be called an open disc 
(of radius $R$).
%  (although the latter is also closed with respect to the
%  topology on $L$ induced by the non-archimedean absolute value $|\;|$). 

\section{Canonical heights  associated to a polynomial map}
\label{globalheightsection}

\subsection{Global canonical heights}

In this subsection, we define the global canonical height function 
associated to
a polynomial $\phi \in K[z]$ of degree $d \geq 2$, and we summarize
some of its main properties.
Our formulation is essentially the same as that in
\cite{Call-Silverman}, although our notation and terminology are slightly different.

If $v \in M_K$, we define the {\it standard local height function} on
$\C_v$ to be the function $h_v : \C_v \to \R$ given by 
\[
h_v(z) = \log^+ | z |_v,
\]
where for $x \in \R$ we set $\log^+ x = \log \max \{ x,1 \}$.

If $S$ is any finite, Galois-stable subset of $\Kbar$, we define the
{\it (absolute logarithmic) height} of $S$ to be the average of the
(properly normalized) local heights of all elements in $S$ (which makes
sense because all but finitely many of these local heights will be
zero).  More precisely, we define
\[
h(S) := \frac{1}{\# S} \sum_{z \in S} \sum_{v \in M_K} \log^+ \| z \|_v.
\]

This definition makes sense, since we have fixed embeddings
of $\Kbar$ into $\C_v$ for each $v$.
However, since we are assuming that $S$ is Galois-stable, the height of $S$ is
in fact independent of our choices of embeddings.

More generally, if $T$ is {\it any} finite subset of $\Kbar$, we define
$h(T)$ to be the height of the smallest Galois-stable subset containing
$T$.  In particular, note that this definition is compatible with the previous one in
the special case where $T$ is already Galois-stable.
If $T = \{ \alpha \}$ consists of a single algebraic number,
the height $h(\alpha)$ we have just defined coincides with the usual
definition of the absolute logarithmic height (see e.g. \cite{SerreLMW}).

\medskip

We now define the (global) canonical height associated to 
$\phi \in K[z]$ by the formula
\[
\hhat_\phi (T) := \lim_{n\to\infty} \frac{1}{d^n} h(\phi^n(T)),
\]
where $\phi^n$ denotes the $n$th iterate of $\phi$.  It follows from
\cite{Call-Silverman} that the limit exists, and in particular that
$\hhat_\phi$ is well-defined.  (This of course uses the fact that $d
\geq 2$).

\medskip

As a function from $\Kbar$ to $\R$, we can characterize $\hhat_\phi$
as the unique function such that 

\begin{itemize}
\item[(1)] There exists a constant $M$ (depending only on $\phi$) such that
  for all $\alpha \in \Kbar$, we have $|\hhat_\phi(\alpha) -
  h(\alpha)| \leq M$.
\item[(2)] For all $\alpha \in \Kbar$, we have $\hhat_\phi(\phi(\alpha)) =
 d\hhat_\phi(\alpha)$.
\end{itemize}

Since $K$ satisfies the Northcott finiteness property, these two
properties imply:
\begin{itemize}
\item[(3)] If $\alpha \in \Kbar$, then $\hhat_\phi(\alpha) \geq 0$, and 
$\hhat_\phi(\alpha) = 0$ if and only if $\alpha$ is a {\it preperiodic point}
for $\phi$, i.e., a point whose orbit $\{ \phi^n (\alpha) : n\geq 0\}$ under
$\phi$ is finite.  
\end{itemize}

\subsection{Canonical local heights}
\label{localheightsection}
There is a 
decomposition of the global canonical height attached to $\phi$
into a sum of canonical local  heights.  In our situation, this is quite
straightforward.\footnote{Note however that one can define canonical
heights attached to more general morphisms, such as a rational
function $\psi$ on $\P^1$, but the canonical local heights are in general
not given by such a simple formula.  This is due to the
fact that as divisors on $\P^1$, we have $\phi^*(\infty) = d(\infty)$,
but $\psi^*(\infty)$ and $d(\infty)$ are merely linearly equivalent.} 
We define $\hhat_{\phi,v} : \C_v \to \R$ to be the function given by
\[
\hhat_{\phi,v}(z) = \lim_{n \to \infty} \frac{1}{d^n} \log^+ | \phi^n(z) |_v.
\]
It is shown in \cite{Call-Silverman} that this limit exists, and that we have
\[
\hhat_\phi (\alpha) = \frac{1}{m}\sum_{v\in M_K} \sum_{i=1}^m n_v 
\hhat_{\phi,v} (\alpha_i),
\]
for all $\alpha\in\Kbar$, where $\alpha_1,\ldots,\alpha_m$ are 
the conjugates of $\alpha$ over $K$.

We can similarly define the local height $\hhat_{\phi,v}(T)$ 
for any finite subset $T$
of $\C_v$.
%\[
% \hhat_{\phi,v}(T) := \frac{1}{\# T} \sum_{z \in T} \hhat_{\phi,v}(z).
% \]

\medskip

If $\phi(z) = a_0 + a_1z + \cdots + a_d z^d  \in K[z]$ with $a_d \neq
0$, and if $v$ is non-archimedean, we let $\alpha_v$ be the
minimum of the $d+1$ numbers $\frac{1}{d-i}v(\frac{a_i}{a_d}), \;
0\leq i < d$ and $\frac{1}{d-1} v(\frac{1}{a_d})$.  (We consider
$v(0)$ to be $+\infty$).

We also define 
\[
c_v(\phi) := |a_d|_v^{\frac{-1}{d-1}}.
\]
See Theorem~\ref{thm1}
for the motivation behind this definition.

The following are some of the properties of the canonical local
heights attached to $\phi$
(see \cite{Call-Silverman} and \cite{Call-Goldstine} for 
proofs):

\medskip

\begin{list1}
\item For each $v$ and each $z \in \C_v$, $\hhat_{\phi,v}(\phi (z)) =
d \hhat_{\phi,v}(z)$.
\item For each $v$, $\hhat_{\phi,v} : \C_v \to \R$ is a continuous
nonnegative function which is zero precisely on the {\it $v$-adic
filled Julia set} $F_v := \{ z \in \C_v \; : \; |\phi^n(z)|_v \not\to \infty
\}$ of $\phi$.
\item For each $v$, the difference $|\hhat_{\phi,v} - h_v|$ is a bounded
function on $\C_v$.
\item For all but finitely many $v \in M_K$, we have
$\hhat_{\phi,v} = h_v$.
\item If $v$ is archimedean, then
$\hhat_{\phi,v}(z) = \log |z|_v - \log c_v(\phi) + o(1)$
as $|z|_v \to \infty$.
\item If $v$ is non-archimedean, then for
$|z|_v$ sufficiently large (depending only on $\phi$), we have
$\hhat_{\phi,v}(z) = \log |z|_v - \log c_v(\phi).$
Specifically, this formula is valid whenever 
$v(z) < \alpha_v$.
\end{list1}

The basic properties of local and global canonical heights provide a
quantitative form of the following ``local-global principle'' for
dynamical systems:

\begin{prop}
An element $\alpha \in \Kbar$ is a preperiodic point for $\phi$ if and
only if every conjugate of $\alpha$ is contained in the $v$-adic filled Julia set $F_v$
of $\phi$ for all $v \in M_K$.
\end{prop}

This principle can also be easily derived directly from the Northcott
finiteness property.

\section{Statement of main results}
\label{resultsection}

Let $L$ be an algebraically closed local field, and let $\phi \in L[z]$ be a polynomial of degree $d\geq 2$.
The {\it filled Julia set} $F = F(\phi)$ of $\phi$ is
the set of points in $L$ which remain bounded under iteration of
$\phi$, i.e.,
\[
F := \{ z \in L \; : \; |\phi^n (z)| \not\to\infty \}.
\]

Our first result concerns the {\it transfinite diameter} of the set
$F$.  The notion of transfinite diameter makes sense in a fairly general
context.  Let $(M,d)$ be a metric space, and let $A$ be a subset of $M$.
We define the $n$th diameter of $A$, denoted $\d_n(A)$, to be the supremum
over all $n$-tuples of points in $A$ of the average pairwise
distance between points
of a given $n$-tuple.  Here ``average'' means geometric mean; in other
words:
\[
\d_n(A) = \sup_{x_1,\ldots,x_n \in A} \prod_{i \neq j} d(x_i,x_j)^
{\frac{1}{n(n-1)}}.
\]

The sequence of nonnegative real numbers
$\{ \d_n(A) \}$ is decreasing, and therefore has a limit as $n$ tends to
infinity.  This is proved, for example, in \cite{Amice}.  
By definition, the {\it transfinite diameter}\footnote{We use the letter $c$ for 
transfinite diameter because $d$ will be the
degree of a polynomial $\phi$, and because over the complex numbers, the
transfinite diameter of a compact set $F$ is the same thing as its
capacity.}
of $A$ is

\[
c(A) = \lim_{n\to\infty} \d_n(A).
\]

By convention, we say that $c(\emptyset) = 0$.

\medskip

In the next section we will compute the transfinite
diameter $c(F)$ of the filled Julia set of a polynomial $\phi \in L[z]$.  
The result is the following:

\begin{theorem}
\label{thm1}
Let $\phi(z) = \sum_{n=0}^d a_n z^n \in L[z]$
be a polynomial of degree $d \geq 2$, and let $F \subset L$ be its filled
Julia set.  Then 
\[
c(F) = |a_d|^{\frac{-1}{d-1}}.
\]
\end{theorem}

In particular, suppose $\phi \in K[z]$ for a global field $K$.  
Then for each $v \in M_K$, one has a 
filled Julia set $F_v = F_v(\phi)
\subset \C_v$ of $\phi$ considered as a polynomial in $\C_v[z]$.
Let $ \c(F_v) = c(F_v)^{n_v}$.
Then as a corollary of Theorem~\ref{thm1}, we obtain 
(using the product formula for $K$):

\begin{cor}
\label{cor1}
If $\phi \in K[z]$ is a polynomial of degree $d \geq 2$, then
\[
\prod_{v\in M_K} \c(F_v) = 1.
\]
\end{cor}

\noindent For notational convenience, we will sometimes write
$c_v(\phi) = |a_d|_v^{-1/(d-1)}$ instead of $c(F_v)$.

As in \cite{Rumely}, the above product formula for the local transfinite
diameters leads to an archimedean equidistribution
theorem for points of small canonical height.  We also obtain, under
additional hypotheses, a family of non-archimedean equidistribution results.
\footnote{We cannot
  directly apply the ideas in \cite{Rumely}, however,
  because the sets $F_v$ are not in general compact.}
In order to state these results, we need to introduce the notion of
an equilibrium measure.  
% Note: Ransford gives a proof of the existence and uniqueness of
% equilibrium measures and includes a discussion of the connection between
% energy integrals and transfinite diameter.

Let $v \in M_K$ be a fixed place of $K$. 
For any compact subset $F \subset \C_v$ with
nonzero transfinite diameter, there exists a unique probability
measure $\mu_F$ which ``minimizes energy'' among all probability measures
supported on $F$.  The energy of a (Borel) probability measure $\mu$ on
$F$ is the double integral
\[
\int_F \int_F -\log |x - y|_v \;d\mu(x) d\mu(y).
\]
The measure $\mu_F$ of minimal energy is
called the {\it equilibrium measure} for $F$, and the energy of $\mu_F$
turns out (not surprisingly) to be $-\log c(F)$.  
(For a more detailed discussion, as well as for proofs of the above 
assertions see  \cite[Chapters 3 and 5 ]{Ransford} for the complex
case and \cite[\S 3.1 and \S 4.1]{Rumelybook} for both 
archimedean and non-archimedean cases).  
% \fixme{Maybe we should omit the reference from \cite{Ransford}}
% No, I don't think that's a good idea.

In particular, if $F$ is compact in $\C_v$ and is the filled Julia set 
for a  polynomial $\phi\in \C_v[z]$, 
then $c(F) = c_v (\phi) > 0$ by Theorem~\ref{thm1}, 
so there
exists a well-defined equilibrium measure $\mu_F$.  

When $\C_v = \C$ and $F$ is the filled Julia set of $\phi$, the support of $\mu_F$
is precisely the topological boundary $\partial F$ of $F$ 
(which coincides with the {\it Julia
set} of $\phi$ -- see Section~\ref{dynamicssection}).
Also, $\mu_F$ is a $\phi$-invariant measure, in the sense that 
\[
\mu_F(\phi^{-1}(B)) = \mu_F(B)
\]
for all Borel sets $B \subseteq \partial F$.
(See \cite[Chapter 6]{Ransford} for these and other facts about
$\mu_F$ in the complex case).

\medskip

Let $M$ be a metric space.
%  which is assumed to be separable and complete.
A sequence  $ \{ \mu_n \}$  of  measures on $M$ is said to 
{\it converge weakly} to a measure $\mu$
%, denoted by $\mu_n\Rightarrow \mu$, 
if
\[
\lim_{n\to \infty}\,\int_{M} \,f d\mu_n = \int_{M}\, f d\mu
\]
for every bounded, continuous function $f : M \to \R$.
% If we only consider measures of compact support, then the word
% *bounded* can be omitted.
For a  general theory of weak convergence of probability 
measures on metric spaces, see \cite{Parthasarathy}.
In this paper, the metric space $M$ will always be a valued 
local field with the metric induced by the absolute value.  

\medskip

For any finite subset $S \subseteq \C_v$,
we denote by $\delta_S$ the probability
measure $\delta_S := \frac{1}{\# S} \sum_{z \in S} \delta_z$, where
$\delta_z$ is the Dirac probability measure supported at the single point
$z \in \C_v$.
Finally, we say a sequence $ \{ S_n \}$ of finite subsets of $\C_v$ is
{\it equidistributed} with respect to the measure $\mu$ if the sequence
$\{ \delta_n \} := \{ \delta_{S_n} \}$ of probability measures converges
weakly to $\mu$.

\medskip

We now introduce a variant of the notion of equidistribution.

First, suppose $F$ is a subset of the valued field $L$, and define a
{\it logarithmic distance function\footnote{{\rm Over $\C$, this is similar 
to the
  notion of a Green's function for $F$, but we do not impose any
  harmonicity condition.}}
 for $F$}
to be a continuous
function $\lambda : L \to \R$ such that $\lambda(z) = 0$ for $z \in F$, 
$\lambda(z) > 0$
for $z \not\in F$, and such that the difference $|\lambda(z) -
\log^+|z| |$ is bounded.

As our main example of such a function, if $\phi \in K[z]$ has degree
$d \geq 2$, then for each place $v \in M_K$ the canonical local height 
function
$\hhat_{\phi,v} : \C_v \to \R$ serves as a logarithmic distance function 
for the
$v$-adic filled Julia set $F_v$ (see Section~\ref{localheightsection}).
If we do not specify a particular set $F$, we define a 
{\it logarithmic distance
function for $L$} to be any nonnegative function $\lambda : L \to \R$
such that the difference $|\lambda(z) - \log^+|z| |$ is bounded.

\medskip

Now let $S = \{ S_n \}$ be a sequence of finite subsets of $L$.
Given a logarithmic distance function $\lambda$ for $F$, we say that the 
sequence
$S$ is $\lambda$-{\it taut} if the average value of $\lambda$ on $S_n$ 
converges to
zero, i.e., if 
\[
\lim_{n \to \infty} \frac{1}{\# S_n} \sum_{z \in S_n} \lambda(z) = 0.
\]

We say that $S$ is a {\it generalized Fekete sequence} for $F$ if $N_n
:= \# S_n \to
\infty$ as $n \to \infty$ and the limit
\[
\lim_{n\to\infty} \prod_{x,y \in S_n, x\neq y}
|x - y|^{\frac{1}{N_n(N_n - 1)}}
\]
exists and equals the transfinite diameter of $F$.

Finally, we say that $S$ is {\it pseudo-equidistributed} with respect to the
pair $(F,\lambda)$ if $S$ is $\lambda$-taut and is a generalized 
Fekete sequence for $F$.

\medskip

With these definitions in mind, we have:

\begin{theorem}
\label{thm2}
Let $K$ be a global field, let $\phi \in K[z]$ be a polynomial of degree at 
least 2, and for each
$v \in M_K$ let $F_v \subset \C_v$ be the $v$-adic filled Julia set of $\phi$.
Let $ \{ S_n \}$ be a sequence of distinct
finite Galois-stable subsets of $\Kbar$
such that $\lim_{n\to \infty} \hhat_\phi (S_n) = 0$.
Then for all $v \in M_K$, the sequence $\{ S_n \}$
is pseudo-equidistributed with respect to the pair $(F_v,\hhat_{\phi,v})$.
\end{theorem}

Over the complex numbers, we will also prove the following result
relating the notions of equidistribution and pseudo-equidistribution:

\begin{prop}
\label{balancedprop}
Let $F \subset \C$ be a compact set with nonzero transfinite diameter,
and let $\lambda$ be a logarithmic distance function for $F$.
Let $S =  \{ S_n \}$ be a sequence of finite subsets of $\C$
which is pseudo-equidistributed with respect to the pair $(F,\lambda)$. 
Then $\{ S_n \}$ is equidistributed with respect to the equilibrium measure on $F$.
\end{prop}

Combining this result with Theorem~\ref{thm2}, we obtain a
generalization of Bilu's equidistribution theorem for sequences of
small points with respect to the usual absolute logarithmic Weil
height \cite{Bilu}.

\medskip

In the non-archimedean case,
if the filled Julia set $F_v \subseteq \C_v$ is compact, 
then we can also prove the following result:

\begin{prop}
\label{nonarch-balancedprop}
Let $L$ be a non-archimedean 
algebraically closed local field and let $\phi\in L[z]$
be a polynomial of degree $d\ge 2$.  Assume that the filled Julia set
$F$ of $\phi$ is compact, and let
$\lambda$ be the canonical local height as defined in 
Section~\ref{localheightsection}. 
Let $S =  \{ S_n \}$ be a sequence of finite subsets of $\C_v$
which is pseudo-equidistributed with respect to the pair $(F,\lambda)$. 
Then $\{ S_n \}$ is equidistributed with respect to the equilibrium 
measure $\mu_F$.
\end{prop}

Combining Theorem~\ref{thm2} and Proposition~\ref{balancedprop}, 
\ref{nonarch-balancedprop}, we
obtain:

\begin{cor}
\label{thm2cor}
Let $K$ and $\phi$ be as in Theorem~\ref{thm2}, and 
let $F \subset \C_v$ be the $v$-adic filled Julia set of $\phi$ with
respect to some fixed  place $v \in M_K$. Assume that $F$ is 
compact with respect to the $v$-adic topology.  
Then if $ \{ S_n \}$ is a sequence of distinct
finite Galois-stable subsets of $\Kbar$
such that $\lim_{n\to \infty} \hhat_\phi (S_n) = 0$,
the sequence $\{ S_n \}$ is equidistributed with respect to $\mu_F$.
\end{cor}

% This generalizes the main result of \cite{Bilu}.

One can apply Theorem~\ref{thm2} and Corollary~\ref{thm2cor} to prove
interesting statements about the canonical height associated to
$\phi$. 
As an example, we present in Section~\ref{Applications}
a generalization to dynamical heights of a result of Schinzel concerning
the usual (logarithmic Weil) height of totally real algebraic numbers.
(See Proposition~\ref{Schinzelresult}). We also present a $p$-adic 
analogue of this result (Theorem~\ref{padicSchinzel}), which
generalizes a theorem of Bombieri and Zannier 
% \cite{Bombieri-Zannier}
for the usual height.

As another application, we show that 
there exists a positive constant $C > 0$ such
that $\hhat_\phi(\alpha) + \hhat_\phi(\sigma
(\alpha)) \geq C$ 
for all but finitely many $\alpha \in \Qbar$ whenever $\sigma$ is
an affine linear transformation of the complex plane which does not
belong to the symmetry group of the complex Julia set of $\phi$
(Theorem~\ref{Bogomolov}). 
This allows us to confirm an interesting special case of a conjecture of S.~Zhang.

% This generalizes a result of S.~Zhang 
% in the special case of lines in affine plane for
% the usual height on $\Qbar$ (see \cite{Zagier} for a discussion).

In Section~\ref{Mandelbrot}, we extend our results to sets which are more
general than the filled Julia sets of polynomial mappings, and as an
application we establish an equidistribution result for the ``adelic
Mandelbrot set''.

\section{Dynamical systems}
\label{dynamicssection}

In this section, we review some terminology and results from
the theory of dynamical systems.  We restrict our attention to 
polynomial dynamics, and for the reader's benefit we give some 
illustrative examples. 

Let $L$ be an algebraically closed local field,
and let $\phi \in L[z]$ be a polynomial.  
In Section~\ref{resultsection}, we defined the filled
Julia set $F$ of $\phi$.  One can give another equivalent definition
of $F$ as follows. 
First of all, it is easy to see that there exists a positive real number $R$
such that if $|z| > R$ then $|\phi(z)| > R$, so that the region 
$\{ z \in L \; : \; |z| > R   \}$ is $\phi$-stable and contained in
the complement of $F$. Let $F_0 := \{ z \in L \; : \; |z| \le R   \}.$
We can then define a sequence $\{ F_n \}$ of  closed subsets
by setting 
$F_{n} = \phi^{-1}(F_{n-1}) = \phi^{-n}(F_0) $ for all $n \ge 1$.
Then
\[
F_0 \supseteq F_1 \supseteq \cdots \supseteq F_n\supseteq \cdots
\supseteq F
\]
and
\[
F = \cap_{n=0}^\infty F_n.
\]
Throughout this paper, we will implicitly choose a decreasing sequence $\{ F_n \}$ 
of   closed subsets as above for the filled Julia set $F$ of $\phi$.
Furthermore, if $L$ is a non-archimedean field then we will always
assume that $R$ is chosen so that
$|\,\phi(z)\,| = |\,a_d z^d\,| > R$ for $|\,z\,| > R$. 

It is clear from the above construction that the filled Julia set is
a closed and bounded subset of $L$. 
If $L = \C$ then $F$ is always compact, but if $L$ is non-archimedean,
then $F$ might or might not be compact.

\medskip

The {\it Fatou set} of $\phi$ is the set of all points in $L$ having a
neighborhood on which the family of iterates $\{ \phi^n \}$ is
equicontinuous (see \cite{Hsia2} for a more detailed discussion of this
definition).  The {\it Julia set} of $\phi$ is defined to be the
complement of the Fatou set. When $L = \C$, the Julia set of $\phi$ is
never empty and is a compact subset of $\C$,
but for non-archimedean
$L$, the Julia set of $\phi$ can be empty or non-empty and non-compact
(see Example~\ref{dynamicsexamples} below). 
Intuitively, the Julia set is the locus
of points which behave ``chaotically'' under iteration of $\phi$.
Moreover, when the Julia set is non-empty, it is the topological 
boundary of the filled 
Julia set (and hence the name of the latter). 
When $L = \C$, this is well-known; when
$L$ is non-archimedean, see \cite[Theorem~5.1]{Benedetto-Component}.

By definition, the Fatou set of $\phi$ is open.  Moreover, 
one can show that the Fatou set, the Julia set, and the filled Julia set of $\phi$ 
are each
{\it completely invariant} under $\phi$, i.e., that each is stable under both
forward and backward iteration of $\phi$.

\medskip

Now we turn to a discussion of periodic and preperiodic points.

\medskip

Let $P \in \Lbar$ be a point of exact period $n$ for $\phi$, i.e., a point
such that $\phi^n(P) = P$ but $\phi^m(P) \neq P$ for all positive
integers $m < n$.  

The {\it multiplier} of $P$ is defined to be $\lambda(P) =
|(\phi^n)'(P)|$.  It follows from the chain rule that $P$ and
$\phi(P)$ have the same multiplier.

We say that $P$ is {\it repelling} (resp. {\it neutral, attractive}) 
if its multiplier satisfies $\lambda(P) > 1$ (resp. $\lambda(P) = 1, \lambda(P) < 1$).

\begin{lemma}
\label{multiplierlemma}
Every attracting periodic point is contained in the Fatou set of
$\phi$, and every repelling periodic point is contained in the Julia
set.  If $L$ is non-archimedean, then every neutral
periodic point is also contained in the Fatou set.  
\end{lemma}

\begin{proof}
See  \cite[Proposition 1.1]{Benedetto-Hyperbolic}.
\end{proof}

For the next result, we introduce another definition.  A point $P \in
\Lbar$ is said to be {\it grand orbit finite} (or an {\em exceptional point})
if the set of all
forward and backward iterates of $P$ under $\phi$ is finite.  Such
points are quite rare, as the following shows.

\begin{lemma}
\label{preperiodiclemma}
If $\phi \in L[z]$ is a polynomial of degree at least 2, then $\phi$
has at most one grand orbit finite point in $\Lbar$.
\end{lemma}

\begin{proof}
See \cite{Beardon} for the case 
$L = \C$ and \cite[Remark. 2.7]{Hsia2} for the non-archimedean case. 
%\helpme{For non-archimedean $L$, this certainly isn't in Beardon!  We
%  need to give a better reference (or a proof).}
\end{proof}

\begin{cor}
\label{preperiodiccor}
$\phi$ has infinitely many distinct preperiodic points in $\Lbar$.
\end{cor}

\begin{proof}
The set $P$ of all preperiodic points of $\phi$ is clearly completely invariant under
$\phi$.  If $P$ is a finite set, then all preperiodic points are
grand orbit finite, and hence there can be only one preperiodic point.  However,
a simple calculation shows this is impossible:   We may assume without
loss of generality that the only preperiodic point of $\phi$ is $0$,
and then we must have $\phi(z) - z  = z^d$, so that $\phi(z) = z^d +
z$.  But then the equation $\phi(z) = 0$ has a nonzero root in $\Lbar$
which is also preperiodic, a contradiction.
\end{proof}

%\helpme{Maybe find a better proof of this lemma, preferably not using the
%  result about grand orbit finite points?  For example, there may be
%  be a simple proof based on counting critical points of $\phi$
%  (although I'm not sure about the inseparable case\ldots)}

We now turn to some illustrative examples of Julia sets and filled Julia sets. Because 
numerous examples of such sets appear
in most books on complex dynamics, we concentrate
below on examples from non-archimedean dynamics.

\begin{example}
\label{dynamicsexamples}
We give examples of polynomials $\phi$ for which the Julia set $J$ of
$\phi$ is (1) empty, (2) nonempty and compact, and (3) non-compact.
\\
\medskip
\begin{itemize}
\item[(1)]
Let $L$ be an algebraically closed non-archimedean local field.  
We say that a polynomial $\phi(z) \in L[z]$ has {\it good reduction}
if $\phi(z) = a_d z^d + a_{d-1}z^{d-1} + \cdots + a_0 \in L[z]$
with $|a_n| \leq 1$ for all $n$ and 
$|a_d| = 1$.  (See \cite{Morton-Silverman} for a more general
definition of good reduction for rational maps). 

It is a simple exercise to verify that if $\phi(z)$ has good
reduction, then the filled Julia set of $\phi$ 
is just the closed unit disc  $\{ |z| \leq 1 \}$ in $L$, and the Julia
set of $\phi$ is empty (something which cannot happen over $\C$!). 

% In the case where ${\mathrm char}(L) = p > 0$ and 
% $\phi(z) = a_0 z + a_1 z^{p} + \cdots + a_d z^{p^d}$ 
% with $|a_0| \le 1$  (no restriction
% on other coefficients), then the filled Julia set
% is not necessary  the unit disc. 

%In this situation, Theorem~\ref{thm1} says that the
%transfinite diameter of the closed unit disc is 1, which is an easy
%computation (see Remark~\ref{capacityofadisk} below).  
%The Julia set of $\phi$ in this case is empty (which is something that
%cannot happen over $\C$!).  
%This example also shows that the filled
%Julia set of a polynomial need not be compact. 
%See \cite{Morton-Silverman,
%Benedetto-Reduction, Hsia1} 
%for the relationship between Julia set and reduction of rational maps. 

% The filled Julia set of a general polynomial $\phi \in L[z]$ can
% be quite complicated, as some of the examples below illustrate.

\item[(2)]
Let $p$ be a prime number.  
Let $k$ be either $\Q_p$ 
or  the field of Laurent series $\Fp ((t))$.
Let $L$ be  a finite extension of ramification index
$e$ and residue degree $f$ over $k$. Let  $|\quad|$ be the absolute value
on $L$ so that $|p| = 1/p$ or $|t| = 1/p$. Let $q$ be the cardinality
of the residue field of $L$, and fix a uniformizer $\pi$ of $L$, so that
$|\pi| = (1/p)^{1/e}$ and $q = p^f$. 

Let  $g(z)\in \Ocal_L [z]$ be a monic 
polynomial of degree $n$ such that
${\bar g}(z) := g(z)  \pmod{\pi}$ is a separable
polynomial of degree $n$ which splits  over $\F_q$.
Take $\phi(z)$ to be the polynomial $ g(z)/\pi$.
Then the filled Julia set and Julia
set of $\phi$ are both equal to a compact subset $J_{\phi}$ of $\Ocal_L$ 
(see \cite{Hsia1} for more detailed description of $J_{\phi}$). 
In particular, if ${\bar g}(z) = z^q - z$ then
$J_{\phi} = \Ocal_L$. 
For example, the Julia set of $g(z) = (z^p - z)/p$ is the compact
subset $\Z_p \subseteq \C_p$.
% It is clear that $J_{\phi}$ is compact 
% since it is a closed subset of the compact set $\Ocal_L$. 

%One can also construct polynomial maps, defined over $\Q_p$, whose
%Julia sets are compact but not contained in any
%finite extension of $\Q_p$.  

%\fixme{It seems correct that the Julia set of the map 
%$\phi(z) = (z (z - 1)^2)/p$ is compact while it's not 
%contained in any finite extension of $\Q_p$ but contained
%in some extension which is also a discrete value field.}

%The formula in Theorem~\ref{thm1}
%says that the transfinite diameter of $J_{\phi}$ is equal to
%$|1/\pi|^{-1/n} = p^{-1/(e\, n)} $. In particular, if 
%$n = q - 1 = p^f - 1$ then $J_{\phi} = \Ocal_L$ and the transfinite
%diameter of $\Ocal_L$ is equal to $ p^{-1/e (p^f - 1)}.$ 
%in agreement with the computation found in
%\cite[Example~4.1.24]{Rumelybook}
% (which relies on a completely different method).

\item[(3)]
Let  $\phi(z) = (z^3 - z^2)/p$.  
If we denote by $v(\cdot)$ the
standard valuation on $\C_p$, so that $v(p) = 1$ and $|z| = p^{-v(z)}$,
then it is easy to see that the filled Julia set $F$ of $\phi$ is contained 
in the unit disc $B_0 := \{
v(z) \geq 0 \}$ and contains the disc $B_\infty := \{ v(z) \geq 1 \}$.
Note that since $F$ contains a disc, it is not compact.  The disc
$B_\infty$ is contained in the Fatou set of $\phi$, but nonetheless
(at least if $p \neq 3$) we claim that the Julia set of $\phi$ is not 
compact.  
To see this, note that by \cite[Example 4.2]{Hsia2}, $\phi$ has both
infinitely many repelling and infinitely many non-repelling periodic
points which are contained in $B_0$.  
(Note that over $\C$, a rational map  of degree at least 2
always has finitely many non-repelling periodic points).
The claim then follows from the following fact (see 
\cite[Proposition 9]{Bezivin} for a proof):  If $P \in \C_p[z]$ is a polynomial
of degree at least 2 whose Julia set is non-empty and compact,  
then all periodic points of $P$ are repelling.

If we set $B_n := \{ v(z) \geq 1 - \frac{1}{2^n} \}$ for
$n=0,1,2,\ldots$, then a simple computation shows that for $n\geq 1$ 
we have $\phi(B_n) \subseteq B_{n-1}$, so that $\phi^n(B_n) \subseteq B_0$.  
Additionally, if we set $C_n := \{ v(z) = 1 - \frac{1}{2^n}  \} $,
then for $z \in \{ 1 - \frac{1}{2^n}< v(z) <  1 - \frac{1}{2^{n+1}}\}
=  B_n\setminus (B_{n+1} \cup C_n)$ ($n \geq 0$), we have
$\phi^{n+1}(z) \not\in B_0$.
In particular, $F$ is contained in the union of $B_\infty$ and
$\cup_{n=0}^\infty C_n$, and the Julia set is contained in 
$\cup_{n=0}^\infty C_n$. 

\end{itemize}
\end{example}

\section{Transfinite diameters}

In this section, we will discuss in more detail the general theory of
transfinite diameters.

\subsection{Fekete's theorem and Chebyshev constants}
\label{Fekete-Chebyshev}

The first result which we want to discuss is a classical theorem of
Fekete (which is the easy half of the ``Fekete-Szego theorem'' -- see
\cite{Rumelybook} for further details).  Since the proof is short and
illustrative, 
%and since we could not find a reference for the exact
%formulation which we give here, 
we include a proof (c.f. \cite[Theorem~5.1.2]{Cantor} and 
\cite[Theorem~6.3.1]{Rumelybook}).

For the statement, we make the following definition.  Let $K$ be a
global field.  An {\it adelic set} with respect to $K$ is a collection
$\{ F_v \}$ of subsets $F_v \subseteq \C_v$ such that $F_v$ is the
closed unit disc in $\C_v$ for almost all $v$.  

We define the transfinite diameter (or capacity) of an adelic set $\FF
= \{ F_v \}$ to be
\[
\c(\FF) := \prod_v \c(F_v).
\]

By Remark~\ref{capacityofadisk} below, almost all terms of this infinite 
product are 1.

%If $K$ is a global field, we define an {\it adelic set}
%with respect to $K$ to be a collection of bounded subsets 
%$F_v \subseteq \C_v$, indexed by the places $v \in M_K$, 
%such that $F_v$ equals the closed unit disc in $\C_v$
%for all but finitely many $v$. For a definition of  
%more general adelic sets of $\P^1$, see \cite{Cantor} and
%\cite{Rumelybook} on adelic sets of  algebraic curves.

%The (global) transfinite diameter of an adelic set $\FF =
%\{ F_v \}$ is  $\prod_{v \in M_K} \c(F_v)$, where $\c(F_v) =
%c(F_v)^{n_v}$.
%This definition makes sense because 
%(see Remark~\ref{capacityofadisk} below) the
%transfinite diameter of the closed unit disc in $\C_v$ equals 1 for
%all $v \in M_K$.

%With this terminology in place, Fekete's result can be stated in the
%following way:

\begin{prop}[Fekete]
\label{Feketetheorem}
Let $K$ be a global field, and let $\FF := \{ F_v \}$ be an adelic set with respect
to $K$ such that $c(\FF) < 1$.  Then there are only finitely many
elements $\alpha \in \Kbar$ such that all Galois conjugates of $\alpha$ are
contained in $F_v$ for all $v \in M_K$.
\end{prop}

\begin{proof}
Assume to the contrary that there are infinitely many such elements
$\alpha_1,\alpha_2,\ldots$, and let $S_n$ be the set of all Galois conjugates
of $\alpha_n$.  By the Northcott finiteness property, the cardinality
$N_n$ of $S_n$ tends to infinity as $n\to\infty$, because the heights
of the $\alpha_n$'s are uniformly bounded by assumption.
Let $\Delta_n := \prod_{x,y \in S_n, x \neq y} (x-y)$ 
be the discriminant of $\alpha_n$, which is a nonzero element of $K$.  
The product formula for $K$ implies that
$\prod_{v \in M_K} \| \Delta_n \|_v = 1$ for all $n$,  so that by the
definition of transfinite diameter we have
\[
1 = \limsup_{n\to\infty} 
\prod_{v \in M_K} \| \Delta_n \|_v^{\frac{1}{N_n(N_n - 1)}}
\leq \prod_{v \in M_K} \c(F_v) < 1,
\]
a contradiction.
\end{proof}

We see that a certain amount of mileage can be gained from exploiting
the interplay between the transfinite diameter and the discriminant.
We will use this observation again in the proof of Theorem~\ref{thm2}.
In order to prove Theorem~\ref{thm1}, however, we will also need some
further properties of the transfinite diameter.
The key result is that when the metric space $M$ in question is a valued
field, the transfinite diameter of a bounded set $A \subset M$ 
coincides with what
is sometimes called the ``Chebyshev constant'' of $A$.  This is made precise in the following lemma,
whose proof can be found in \cite[Lemma 5.4.2]{Amice} (see also
\cite[Theorem 16.2.1]{Hille}, which is stated for compact subsets of
$\C$, but whose proof can be easily adapted to the present situation). 

\begin{lemma}
\label{Chebyshev}
Let $B$ be a bounded subset of the valued field $L$, and denote by
$P_n$ the set of monic polynomials of degree $n$ with coefficients in $L$.
If $P \in P_n$, define $ \| P \|_B = \sup_{x \in B} |P(x)|$, and let
\[
S_n(B) = \inf_{P \in P_n} \| P \|_B, \; \; s_n(B) = (S_n(B))^{1/n}.
\]
Then $s_n(B) \to c(B)$ as $n \to\infty$.
\end{lemma}

The proof of the next lemma is adapted from the exercises for 
\cite[Section 16.2]{Hille}.

\begin{lemma}
\label{Hille}
Let $L$ be an algebraically closed valued field and let
$A$ be a bounded subset of $L$. Let $\phi \in L[z]$ be a 
polynomial of degree $d \geq 1$ with
leading coefficient $a_d$, and set
$A^* := \{ z \in L \; : \; |\phi(z)| \in A \}$.
Then $c(A^*) = \left(c(A)/|a_d|\right)^{\frac{1}{d}}$.
\end{lemma}

\begin{proof}
It is an elementary exercise to show that the assertion is true 
in the case $d = 1$. For the general case,  
a simple change of  variables shows that
we only need to treat the case where $\phi$ is monic, in which case we
want to show that $c(A^*) = c(A)^{\frac{1}{d}}$.

Let $T_n \in P_n$ be arbitrary.  Then $T_n \circ \phi \in P_{nd}$ and
$\| T_n \circ \phi \|_{A^*} \leq \| T_n \|_A.$  It therefore follows
from Lemma~\ref{Chebyshev} that $c(A^*) \leq c(A)^{\frac{1}{d}}$.

For the other direction, we again use Lemma~\ref{Chebyshev}.  Let $n
\geq 1$ be an integer, let $\epsilon > 0$ be arbitrary, and
choose $T_n^* \in P_n$ such that $\| T_n^* \|_{A^*} < S_n(A^*) + \epsilon$.
For each $w \in L$, let $z_1(w),\ldots,z_d(w)$ denote the roots of
$\phi(z)=w$ (with multiplicities), and set
$\phi_n(w) := \prod_{j=1}^d T_n^*(z_j(w))$.  Viewing $\phi_n(w)$
as the resultant of $T_n^*(z)$ and $\phi(z) - w$, it is clear from
Sylvester's determinant that $\phi_n$ is a
polynomial of degree $n$ in $w$ with leading coefficient $\pm 1$.
Therefore 
\[
S_n(A) \leq \| \phi_n \|_{A} \leq \prod_{j=1}^d \| T_n^* \|_{A^*}
= \| T_n^* \|_{A^*}^d < (S_n(A^*) + \epsilon)^d.
\]
Since $\epsilon$ was arbitrary, it follows that $s_n(A) \leq s_n(A^*)^d$
for all $n\geq 1$.  Therefore $c(A^*) \geq c(A)^{\frac{1}{d}}$ as desired.
\end{proof}

% Old version of the ``Amice'' lemma:
% \begin{lemma}
% \label{Amice}
% Let $L$ be an algebraically closed valued field, let
% $\phi \in L[z]$ be a monic polynomial of degree $d \geq 1$, let $M$ be a
% nonnegative real number, and let
% $B := \{ z \in L \; | \; |\phi(z)| \leq M \}$.
% Then $c(B) = M^{\frac{1}{d}}$.
% \end{lemma}

\begin{remark}
\label{capacityofadisk}

\begin{itemize}
\item[(1)]
Suppose $L$ is an algebraically closed local field.
If $M$ is a positive real number such that
$|\alpha| = M$ for some $\alpha \in L$ (in which case we say that $M$
is in the {\it value group} of $L$), then the disc $D_M := \{ z \in L \; : \;
|z| \leq M \}$ has transfinite diameter $M$.  To see this, note that
by making the substitution $z \mapsto z/\alpha$, we may reduce to the 
statement
that if $D = D_1$ is the unit disc in $L$, then $c(D) = 1$.  
Since the polynomial $z^n$ has sup-norm equal to 1 on $D$ for all $n$,
it follows from Lemma~\ref{Chebyshev} that $c(D) \leq 1$.  On the
other hand, for any $n$ not divisible by the residue characteristic of $L$, 
consider the
$n$ distinct roots $\zeta_1,\ldots,\zeta_n$ of the separable
polynomial $f_n(z) := z^n - 1$.  We find that 
\[
\prod_{i \neq j} |\zeta_i - \zeta_j| = \prod_i |f_n'(\zeta_i)| = |n^n|
\geq 1,
\]
so that $c(D) \geq 1$ by the definition of transfinite diameter.  Therefore 
$c(D) = 1$ as claimed.

\item[(2)]
\label{PLremark}
Let $\phi \in L[z]$ be a polynomial of degree $d \ge 1$ with leading 
coefficient $a_d$. 
Applying Lemma~\ref{Hille} and (1), we see that the transfinite diameter of the region
$B := \{ z \in L \; : \; |\phi(z)| \leq M \}$ is equal to 
$(\frac{M}{|a_d|})^{\frac{1}{d}}$.
We will refer to such regions as {\em lemniscates}.
% Or maybe as {\em PL-domains}.  (The ``PL'' stands for ``polynomial lemniscate''.)
\end{itemize}
\end{remark}

\subsection{Proof of Theorem \ref{thm1}}
\label{thm1proofsection}

The conclusion of Theorem~\ref{thm1} in the case $L = \C$ is well-known 
(see e.g.
\cite{Ransford}).  However, the proof uses the following result
(see \cite[Theorem 5.1.3]{Ransford} or \cite[Theorem 16.2.2(iii)]{Hille} for a proof), 
which does not generalize to the more general fields $L$ which we consider:

\begin{lemma}
\label{intersectionlemma}
Let $F_1,F_2,\ldots$ be closed and bounded subsets of $\C$ with 
$F_1 \supseteq F_2 \supseteq \cdots$, 
and let $F = \cap_{n=1}^\infty F_n$.  Then 
\[
c(F) = \lim_{n\to\infty} c(F_n).
\]
\end{lemma}

\begin{remark}
The analogue of this lemma is false in general over a non-archimedean valued
field $L$ (such as $\C_p$ for any prime $p$) which is not locally
compact.  
(Of course, one always has an
inequality $c(F) \leq \lim_{n\to\infty} c(F_n)$.)

For example, if $L=\C_p$, we can find an infinite sequence of disjoint
closed discs $B_1,B_2,\ldots$ of the form $B_i = \{ z \in \C_p \; :
\; |z - \alpha_i| \leq \frac{1}{p} \}$ for some sequence $\alpha_i \in \C_p$ with each
$|\alpha_i| = 1$.  Then if we set $F_n := \cup_{i=n}^\infty B_i$ for
$n\geq 1$, then each $F_n$ is closed and bounded and $F_1 \supseteq
F_2 \supseteq \cdots$.
However, we have $c(F_n) \geq c(B_n) = \frac{1}{p}$ for all $n$ (see
Remark~\ref{capacityofadisk}), whereas the transfinite diameter
of $\cap_{n=1}^\infty F_n = \emptyset$ is zero.
\end{remark}

\begin{remark}
\label{nesteddiskremark}
The analogue of Lemma~\ref{intersectionlemma} over $\C_p$ is still
false even if we assume that each set $F_n$ is a closed disc, since
one can construct a sequence $F_1 \supseteq F_2 \supseteq \cdots$
of closed discs of radius $r_n$ in $\C_p$ such that
$\liminf_{n\to\infty} r_n > 0$
but $\cap_{n=1}^\infty F_n = \emptyset$.
\end{remark}

% To do this, put an equivalence relation on $\C_p$ whereby $x \sim y$
% iff $|x-y|_p \leq 1$.  Then since $\C_p$ has a countable dense
% subset (e.g. $\Qbar$), it follows that there are a countable number
% of equivalence classes.  Enumerate these classes as $E_2,E_3,E_4,\ldots$
% Then starting with $F_1 = D(a_1,r_1)$, where $a_1 = 0$ and $r_1 =
% 2$, if we have chosen $a_n$ and $r_n$, we can choose any $r_{n+1}$
% such that $1 < r_{n+1} < r_n$ and then we can choose $a_{n+1}$ so
% that $F_{n+1} := D(a_{n+1},r_{n+1})$ and $E_{n+1}$ are disjoint.  It
% is easy to see that $\{ F_n \}$ has the desired properties.
%\fixme{I still don't quite understand this example.}
%\helpme{Should I explain this more?  I think it's correct.}

We now prove Theorem \ref{thm1}.  

\begin{proof}

Define $c(\phi) = |a_d|^{\frac{-1}{d-1}}$, where $a_d$ is the leading
coefficient of $\phi \in L[z]$.
The proof of the
inequality $c(F) \leq c(\phi)$
is based on the argument found in \cite[Theorem 6.5.1]{Ransford}.
We give two proofs of the reverse inequality, one global and the
other purely local.  The global proof is in some sense simpler, but it
only applies under additional hypotheses.

%Let $a_d$ be the leading coefficient of $\phi$, and choose $\alpha \in L$ such that
%$\alpha^{d-1} = a_d$, and let $m(z) = \alpha z$.  If we let
%$\tilde{\phi} = m \circ \phi \circ m^{-1}$, we see that
%$\tilde{\phi}$ is a monic polynomial of degree $d$ whose filled
%Julia set is just $\tilde{F} = m(F)$.  It therefore suffices to
%prove that if $\phi(z)$ is {\it monic}, then $c(F) = 1$.  Having
%reduced to this special case, we proceed as follows.

Recall that $\phi(z) \in L[z]$, where $L$ is either $\C$ or an
algebraically closed ultrametric field, that the degree
of $\phi$ is $d \geq 2$, and that $F$ is the filled Julia set of
$\phi$. Fix  a sequence $\{ F_n \}$ of  closed subsets
as in Section~\ref{dynamicssection} such that  
$F_{n} = \phi^{-1}(F_{n-1}) = \phi^{-n}(F_0) $ for all $n \ge 1$, so that
\[
F_0 \supseteq F_1 \supseteq \cdots \supseteq F_n\supseteq \cdots
\supseteq F
\]
and
\[
F = \cap_{n=0}^\infty F_n.
\]
By Lemma~\ref{Hille}, we have
\[
c(F) \leq \lim_{n\to\infty} c(F_n) = \lim_{n\to \infty}
 c(F_0)^{\frac{1}{d^n}} 
|a_d|^{-(\frac{1}{d} + \cdots + \frac{1}{d^n})} = |a_d|^{\frac{-1}{d - 1}},
\]
which gives us the upper bound $c(F) \leq c(\phi)$.

\medskip

We now present two proofs of the fact that $c(F) = c(\phi)$.  

\medskip

For the global proof, we make the additional assumption that
$\phi \in K[z]$, where $K$ is a global field, and that 
$L = \C_v$ for some $v \in M_K$.
%We may again assume without loss of generality that $\phi$ is monic,
%by enlarging $K$ if necessary to contain an element $\alpha$ with
%$\alpha^{d-1} = a_d$.
We know already that $c(F) \leq c(\phi)$ for all $v \in M_K$, and the idea
is to prove the reverse inequality for all $v\in M_K$ at once.
To do this, we first recall that by Corollary~\ref{preperiodiccor}, 
$\phi$ has infinitely many 
distinct preperiodic points in $\Kbar$, and by definition, 
a preperiodic point of $\phi$ lies in $F_v$ for all places $v$.  

Let $\FF$ be the adelic set corresponding to the collection $\{ F_v \}$.  
If $c(F_v) < c_v(\phi)$ for some $v \in M_K$, then $\c(\FF) < \prod_{v
\in M_K} \c_v(\phi) = 1$ and the
preperiodic points furnish a contradiction to Fekete's theorem (Theorem~\ref{Feketetheorem}).
Therefore $c(F_v) = c_v(\phi)$ for all $v$, as desired.

\medskip

For the purely local proof, we return to the general 
setup where $\phi \in L[z]$ is 
a  polynomial of degree at least 2, and we need to show that $c(F) \geq 
c(\phi)$.  
Note that if $L = \C$, then we can already conclude that $c(F)=c(\phi)$ by 
Lemma~\ref{intersectionlemma}.  So without loss of generality, we may
assume that the field $L$ is non-archimedean.

Since $F$ is completely invariant under $\phi$, we have 
$F = \phi^{-1}(F) = \{ z \in L \; : \; \phi(z) \in F \}$.  By Lemma~\ref{Hille},
 we
have $c(F)^d = c(\phi)^{d - 1} c(F)$, so it suffices to show that $c(F) > 0$.
%  If $F$
%contains a disk (of positive radius), then it is clear from
%Remark~\ref{capacityofadisk} that $c(F) > 0$.  Also, we know from
%Lemmas~\ref{disksinF} and \ref{multiplierlemma} that if $F$ contains a periodic point which is
%either attracting or indifferent, then $F$ contains an entire dis5k
%around that point.
If there is a non-repelling periodic point $Q$
in $F$, then $Q$ is in the Fatou set by Lemma~\ref{multiplierlemma}. 
As the Fatou set is open
and completely invariant under $\phi$, it follows that 
$F$ contains an entire disc around $Q$. It is then clear from  
Remark~\ref{capacityofadisk} (1) that $c(F) > 0$. 
Therefore we may reduce to the case where all periodic points of $F$
are repelling.

For each $n\geq 1$, let ${\rm Per}_n$ denote the set of all points of
period dividing $n$ (i.e., the set of all roots of $\phi^n(z) - z$)
and let $N_n$ be its cardinality.
By the assumption that all periodic points are repelling, together
with the fact that the absolute value on $L$ is non-archimedean, we see that 
$|(\phi^n)'(P) - 1 | = |(\phi^n)'(P)  | > 1$ for all $P \in {\rm Per}_n$. 
In particular, the polynomial $\phi^n(z) - z$ is separable, so that
$N_n = d^n$.
Letting $\d({\rm Per}_n) := \prod_{x,y \in {\rm Per}_n, x \neq y} 
|x - y|^{1/ N_n (N_n - 1)}$,
we find (using a standard algebraic identity) that
\[
\d({\rm Per}_n) = 
|a_d|^{\frac{-1}{d-1}} \prod_{P \in {\rm Per}_n}  
|(\phi^n)'(P) - 1|^{1/ N_n (N_n - 1)}
> c(\phi)  
\]
for all $n$, so that $c(F) \geq c(\phi)$, as desired.
\end{proof}

% In the following we apply Theorem~\ref{thm1} to 
% the filled Julia sets computed in Example~\ref{dynamicsexamples}.

\begin{example}
\label{computetrasfinitediameters}

% \begin{list1}
% \item For the case where $\phi(z) = a_d z^d + a_{d-1}z^{d-1} + \cdots + a_0$
% with $|a_n| \leq 1$ for all $n$ and 
% $|a_d| = 1$, in this situation, Theorem~\ref{thm1} says that the
% transfinite diameter of the closed unit disc is 1 (c.f. 
% Remark~\ref{capacityofadisk}(1)).

Let $\phi(z)$ be as given in Example~\ref{dynamicsexamples}(2).
Then, the formula in Theorem~\ref{thm1}
says that the transfinite diameter of $J_{\phi}$ is equal to
$|1/\pi|^{-1/(n-1)} = p^{-1/(e\, (n-1))} $. In particular, if 
$n = q = p^f $ then $J_{\phi} = \Ocal_L$ and the transfinite
diameter of $\Ocal_L$ is equal to $ p^{-1/e (p^f - 1)}$,
in agreement with the computation found in
\cite[Example~4.1.24]{Rumelybook}.
%(which relies on a completely different method).
As a concrete special case, the transfinite diameter of $\Z_p
\subseteq \C_p$ is $p^{\frac{-1}{p-1}}$.
% \end{list1} 
\end{example}

\section{Equidistribution and pseudo-equidistribution}

\subsection{Proof of Theorem~\ref{thm2}}

We begin by proving a couple of important technical lemmas concerning the
relationship between heights and transfinite diameters.
First, we make some definitions.  
Suppose $L$ is a valued field and that $S = \{ S_n \}$ and $S' = \{ S_n' \}$ 
are sequences of finite subsets of $L$.  
Let $N_n$ (resp. $N_n'$) denote the cardinality of $S_n$
(resp. $S_n'$).

We say that $S'$ is a {\it full subsystem} of $S$ if:
\begin{itemize}
\item $S_n' \subseteq S_n$ for all $n$
\item $\lim_{n \to \infty} N_n/N_n' = 1$.
\end{itemize}

If $T$ is any finite subset of $L$, with cardinality $N$, we define 
\[
\d(T) := \prod_{x,y \in T, x\neq y} |x - y|^{\frac{1}{N(N - 1)}}.
\]

With this terminology, we have:

\begin{lemma}
\label{technicallemma}
Let $L$ be a valued field, let $S = \{ S_n \}$ and $S' = \{ S_n' \}$ be
sequences of finite subsets of $L$, and assume that $S'$ a full subsystem of $S$.  
Let $N_n$ (resp. $N_n'$) denote the cardinality of $S_n$
(resp. $S_n'$), and assume that $\lim_{n\to\infty} N_n = \infty$.
Furthermore, let $\lambda$ be a logarithmic distance function for $L$,
and assume that $S$ is $\lambda$-taut.
% i.e., that $\frac{1}{N_n} \sum_{z \in S_n} \lambda(z) \to 0$ as $n\to\infty$.
Then 
\[
\limsup \d(S_n) \leq \limsup \d(S_n').
\]
\end{lemma}

\begin{proof}
Set $T_n = S_n \setminus S_n'$, then 

%\[
%\d(S_n) = \frac{1}{N_n(N_n - 1)} \sum_{x,y \in S_n', x \neq y} 
%\log |x-y| + 
%\frac{2}{N_n(N_n - 1)} \sum_{x \in S_n,y \in T_n, x \neq y} 
%\log |x-y|.
%\]
\[
\d(S_n) = \frac{1}{N_n(N_n - 1)} \sum_{x,y \in S_n', x \neq y} 
\log |x-y| + \psi_n , 
\]
where 
\[
\psi_n  := \frac{2}{N_n(N_n - 1)} \sum_{x \in S_n' ,y \in T_n} 
\log |x-y| + 
\frac{1}{N_n(N_n - 1)} \sum_{x, y \in T_n, x \neq y} 
\log |x-y|. 
\]

Since $\lim_{n \to \infty} N_n/N_n' = 1$, we find 
using Remark~\ref{DoubleSequenceRemark} below that
\[
\limsup \d(S_n) \le \limsup \d(S_n') + \limsup \psi_n.
\]

Observe that $\;\log|x - y| \leq \log^+ |x|+ \log^+ |y| + \log 2$ 
for all $x, y \in L$, since by the triangle inequality we have
\[ 
|x - y| \leq 2 \max\{|x|,|y|\} \leq 2 \max\{1, |x|\}\max\{1, |y|\}. 
\]

By assumption, there exists a constant $C>0$ such that $\log^+|z| \leq
\lambda(z) + C$.
Therefore $\log |x-y| \leq \lambda(x) + \lambda(y) + C'$,
where $C' = 2C + \log 2$.

If we set $R_n := \# T_n$, so that $\lim R_n/N_n = 0$ by assumption, then
it follows that
\begin{equation*}
\begin{aligned}
\psi_n  &\leq \frac{1}{N_n(N_n - 1)}\left\{ 
2 \sum_{\substack{x \in S_n' \\ y \in T_n}} 
(\lambda(x) + \lambda(y) + C') + 
\sum_{\substack{x, y \in T_n \\ x \neq y}} (\lambda(x) + \lambda(y) + C')\right\}  \\ 
%   & =  \frac{2}{N_n} \left\{\frac{R_n}{N_n - 1}\sum_{z \in S_n'} \lambda(z) 
% + \sum_{z\in T}\, \lambda(z) + R_n C'\right\}
&=  \frac{1}{N_n(N_n - 1)} \left\{
2R_n \left( \sum_{z \in S_n'} \lambda(z) \right)
 + 2(N_n' + R_n - 1) \left( \sum_{z\in T_n}\, \lambda(z) \right)
\right\} \\
&+  \frac{1}{N_n(N_n - 1)} \left\{ R_n(2N_n' + R_n - 1) C' \right\} . 
\end{aligned}
\end{equation*}

Therefore $\limsup \psi_n \leq 0$ by our assumptions, which yields the
desired inequality.
\end{proof}

The following key lemma makes precise the intuition that if a set $T$ of points in $\C_v$
is close to the filled Julia set $F_v$ of $\phi$ (in the sense that
the canonical local height $\hhat_{\phi,v}(T)$ is small\footnote{
We remind the reader that by definition, 
$\hhat_{\phi,v}(T) := \frac{1}{\# T} \sum_{z \in T}
\hhat_{\phi,v}(z).$
Note that one should not take the above intuition too literally, however,
since a sequence of points $z_n \in \C_v$ for $v$ non-archimedean can satisfy
$\hhat_{\phi,v}(z_n) \to 0$ without ${\rm dist}(z_n,F_v)$ converging to zero
(see Example~\ref{localheightexample} below).}), 
then the average pairwise distance between the points of $T$ can exceed the transfinite
diameter of $F_v$ only by a small amount.

\begin{lemma}
\label{keylemma}
Let $v\in M_K$, and suppose that $S_n$ is a sequence of
finite subsets of $\C_v$ with the property that
$N_n = \# S_n \to \infty$ and  $\hhat_{\phi,v}(S_n) \to 0$.
Then
\[
\limsup_{n\to\infty} \d(S_n) \leq c_v(\phi).
% Recall that \d(S_n)
% \prod_{x,y \in S_n, x\neq y}
% |x - y|_v^{\frac{1}{N_n(N_n - 1)}} 
\]
\end{lemma}

\begin{proof}

Set $\lambda(z) := \hhat_{\phi,v}(z)$ and put $L := \C_v$.  By property (3) above of the
local canonical height, $\lambda$ is a logarithmic distance function
for $L$, so that there exists a constant $C>0$ such that 
$|\lambda(z) - \log^+|z|| \leq C$ for all $z \in L$.
In particular (with the convention that $\log 0 = -\infty$), we have 
$\log |z| \leq \lambda(z) + C$ for all $z \in L$.

By property (1) of the local canonical height, we also have 
$\lambda(\phi(z)) = d\lambda(z)$ for all $z \in L$.  Therefore, for
all $m\geq 1$ and all $z\in L$ we have
\begin{equation}
\label{heightinequality}
\frac{1}{d^m}\log |\phi^m(z)| \leq \lambda(z) + \frac{C}{d^m}.
\end{equation}

Fix a rational number $\beta > 0$, and define 
$S_n' := \{ z \in S_n \; : \; \lambda(z) \leq \beta / 2 \}$.
Then $S' := \{ S_n' \}$ is a full subsystem of $S$, since (setting
$T_n = S_n \setminus S_n'$ and $R_n = \# T_n$) if $\limsup R_n / N_n = \delta >
0$, then we would have $\hhat_{\phi,v}(S_n) \geq \delta\beta > 0$ for
infinitely many $n$, contradicting the fact that $\hhat_{\phi,v}(S_n) \to 0$.

By Lemma~\ref{technicallemma}, it suffices to prove that
$\limsup_{n\to\infty} \d(S_n') \leq c_v(\phi)$.

\smallskip

To see this, for each integer $m\geq 1$ we define
\[
U_{m,\beta} := \{ z \in L \; : \; \frac{1}{d^m} \log |\phi^m(z)|
\leq \beta   \}.
\]

Then by (\ref{heightinequality}) above, for any $m$ such that
$\frac{C}{d^m} \leq \beta/2$ we have $S_n' \subset U_{m,\beta}$.
By the definition of transfinite diameter, it follows that for $m$
sufficiently large we have
\begin{equation}
\label{Umbetainequality}
\limsup_{n\to\infty} \d(S_n') \leq c(U_{m,\beta}).
\end{equation}

On the other hand, note that $\phi^m(z)$ has degree $d^m$ and leading
coefficient $a(m) := a_d^{1 + d + d^2 + \cdots + d^{m-1}}$.
% Could also say $a(m) := a_d^{\frac{d^m - 1}{d-1}}$.
Since
\[
|a(m)|^{-d^{-m}} = \left( |a_d|^{\frac{d^m-1}{d-1}} \right)^{-d^{-m}}
= \left( |a_d|^{\frac{-1}{d-1}} \right)^{1 - d^{-m}} = c_v(\phi)^{1-d^{-m}},
\]
it follows from Remark~\ref{capacityofadisk}(2) that
\[
c(U_{m,\beta}) = \frac{e^\beta}{|a(m)|^{\frac{1}{d^m}}} 
= e^{\beta} c_v(\phi)^{1 - d^{-m}}.
\]

Letting $m\to\infty$ in (\ref{Umbetainequality}), we find that
$\limsup_{n\to\infty} \d(S_n') \leq c_v(\phi)e^\beta.$
Now let $\beta \to 0$, and we obtain
$\limsup_{n\to\infty} \d(S_n') \leq c_v(\phi)$ as desired.

\end{proof}

We now turn to the proof of Theorem~\ref{thm2}.

\begin{remark}
\label{DoubleSequenceRemark}
Before we give the proof of Theorem~\ref{thm2}, 
note that if $(a_n^{(j)})_{j=1,\ldots,k, n
\geq 1}$ is a doubly indexed sequence of real numbers, then
\[
\limsup_{n\to\infty} ( \sum_{j=1}^k a_n^{(j)} ) \leq
\sum_{j=1}^k \limsup_{n\to\infty} a_n^{(j)}.
\]
This is an easy consequence of Fatou's lemma, applied to an appropriate step function.
However, this inequality can fail if we replace the finite sum by an infinite sum.
\end{remark}

\begin{proof}
Let $N_n$ be the cardinality of $S_n$.  Since $K$ has the Northcott
finiteness property and the difference $|\hhat_\phi (\alpha) - h(\alpha)|$ is bounded
on $\Kbar$, it follows from the assumption $\hhat_{\phi}(S_n) \to 0$ 
that $N_n \to \infty$ as $n\to\infty$.
By hypothesis, we have $\sum_{v \in M_K} \hhat_{\phi,v}(S_n) \to 0$.
Since each of the functions $\hhat_{\phi,v}$ takes nonnegative values,
this implies that
$\frac{1}{N_n} \sum_{z \in S_n} \hhat_{\phi,v}(z) \to 0$
for all $v \in M_K$, i.e., that $S$ is $\hhat_{\phi,v}$-taut for all $v$.

\smallskip

We need to show that $S$ is a generalized Fekete sequence for $F_v$
for all $v \in M_K$.  
First, we note that  $\log \c_v(\phi)= 0 $ and 
$\hhat_{\phi,v}(z) = h_v(z)$ for all but finitely many $v\in M_K$. Let 
${\mathcal T}$ be a finite subset of $M_K$ containing the archimedean 
places and all places $v$ where either $\log \c_v(\phi) \ne 0 $ or 
$\hhat_{\phi,v}(z) \ne h_v(z)$.

%To do this, we
%proceed as in the proof of Theorem~\ref{Feketetheorem}.
% The idea is to now combine the product formula for the local transfinite
% diameters with the fact that, since each $S_n$ is Galois-stable, the product
% $ \Delta_n := \prod_{x,y \in S_n, x\neq y} (x - y)$ is a nonzero $K$-rational
% number.

Let $ \Delta_n := \prod_{x,y \in S_n, x\neq y} (x - y)$, a
nonzero element of $K$. 
For $v\in M_K\setminus {\mathcal T}$, we have 
\[
\log |x - y|_v \le \log^+ |x|_v + \log^+ |y|_v = h_v(x) + h_v(y)
\]
as $v$ is a non-archimedean 
valuation.   Therefore, we can bound $\log\|\Delta_n\|_v$ in terms of
local heights of $S_n$ as follows:
\begin{align*}
\frac{1}{ N_n( N_n - 1)} \log 
\| \Delta_n\|_v & \leq 
\frac{n_v}{N_n(N_n-1)} \sum_{x,y\in S_n, x\ne y} \{ h_v(x) +h_v(y)\} 
\\
  & =  \frac{2 n_v}{ N_n}  \sum_{x\in S_n} h_v(x). 
\end{align*}
% where $h_v(S_n) = \sum_{x\in S_n} h_v(x).$
Hence,
\begin{align*}
\frac{1}{ N_n( N_n - 1)} \sum_{v \notin {\mathcal T}} \log
\|\Delta_n \|_v 
& \leq  \frac{2}{N_n} \sum_{v\notin {\mathcal T}} n_v
\sum_{x\in S_n} h_v(x)
\\
   & \leq  2 \hhat_{\phi}(S_n),
\end{align*}
where the last inequality comes from the 
%choice of ${\mathcal T}$ and
non-negativity of local heights. 

Using Lemma~\ref{keylemma}, we have
\[
\limsup_{n\to\infty} \frac{1}{ N_n( N_n - 1)} \log \| \Delta_n \|_v 
\leq \log \c_v(\phi)
\]
for all $v \in M_K.$ 
By assumption, $\lim_{n\to \infty} \hhat_{\phi}(S_n) = 0$, and 
since $\Delta_n \neq 0$, the product formula gives
$\frac{1}{ N_n( N_n - 1)} \sum_{v \in M_K} \log \| \Delta_n \|_v = 0$
for all $n$. Also note that  by the product formula we have 
$\sum_{v\in {\mathcal T}} \log \c_v(\phi) = \sum_{v\in M_K}\log \c_v(\phi) = 0.$
Therefore,  
\begin{equation}
\label{eq:limsupeqn}
\begin{aligned}
0 & = \limsup_{n\to\infty}  \frac{1}{ N_n( N_n - 1)}
\sum_{v \in M_K} \log \|\Delta_n \|_v \\
&\leq\sum_{v\in {\mathcal T}} 
\limsup_{n\to\infty} \frac{1}{ N_n( N_n - 1)} \log \|\Delta_n \|_v 
  + \limsup_{n\to\infty} 
\frac{1}{ N_n( N_n - 1)} \sum_{v\notin {\mathcal T}} 
\log \|\Delta_n\|_v  \\
  & \leq \sum_{v \in {\mathcal T}} \log \c_v(\phi) + 2 \lim_{n\to
    \infty} \hhat_{\phi}(S_n) = 0.
\end{aligned}
\end{equation}
It follows that $\limsup_{n\to\infty}
\frac{1}{ N_n( N_n - 1)}  \log \|\Delta_n \|_v  = \log\c_v(\phi)$ for
all $v\in {\mathcal T}$. As ${\mathcal T}$ can be chosen to contain
any place $v\in M_K$,  we have shown that
\[
\limsup_{n\to\infty}
\frac{1}{ N_n( N_n - 1)}  \log \|\Delta_n \|_v  = \log\c_v(\phi)
\]
for all $v\in M_K$.  

To finish the proof, let $w\in M_K$ be any fixed place. Then
\[
\log \|\Delta_n \|_w = 
- \sum_{v\in M_K, v\ne w} \log \|\Delta_n
\|_v .
\]

Let ${\mathcal T}$ be a finite set of places containing $w$.  Then by
the same argument as in (\ref{eq:limsupeqn}), we have
\[
\limsup_{n\to\infty}  \frac{1}{ N_n( N_n - 1)}
\sum_{v \in M_K, v \neq w} \log \|\Delta_n \|_v 
\leq \sum_{v \in {\mathcal T}, v \neq w} \log \c_v(\phi)
= -\log \c_w(\phi).
\]

It follows that
\begin{align*}
\liminf_{n \to \infty} \frac{1}{ N_n( N_n - 1)}  \log \|\Delta_n
\|_w & = - \limsup_{n\to \infty}\,\frac{1}{ N_n( N_n - 1)} 
\sum_{v\in M_K, v\ne w} \log \|\Delta_n\|_v \\
& \ge \log\c_w (\phi). 
\end{align*}

% Previous argument:
% \begin{align*}
% \liminf_{n \to \infty} \frac{1}{ N_n( N_n - 1)}  \log \|\Delta_n
% \|_w & = - \limsup_{n\to \infty}\,\frac{1}{ N_n( N_n - 1)} 
%  \sum_{v\in M_K, v\ne v_0} \log \|\Delta_n\|_v \\
%   & \ge - \sum_{v\in M_K, v\ne w} \limsup_{n\to \infty}\,
% \frac{1}{N_n( N_n - 1)} \log \|\Delta_n\|_v \\ 
%   & = - \sum_{v\in M_K, v\ne w} \log\c_v(\phi) =
%   \log\c_w (\phi). 
% \end{align*}
As $w$ is arbitrary, we therefore have 
\[
\lim_{n\to\infty} \frac{1}{ N_n( N_n - 1)} \log \| \Delta_n \|_v = \log\c_v(\phi)
\]
for all $v \in M_K$.
Thus $S$ is pseudo-equidistributed with respect to
$(F_v,\hhat_{\phi,v})$ as desired.
% In particular, $S$ is a generalized Fekete sequence for $F$, as was to 
% be shown.
\end{proof}

% Intuitively, one might think that if a sequence $S = \{S_n\}$ of
% finite subsets is pseudo-equidistributed with respect to 
% $(F_v, \hhat_{\phi,v})$, then they should spread over $F_v$ as 
% possible as they can. However, the following simple example
% indicates that in general pseudo-equidistribution may not be 
% \lq\lq equidistributed \rq\rq . 
% \fixme{My comment might be wrong and probably need to be 
% corrected.}

\begin{example}
Let $p$ be a prime number. 
Assume that $K$ is a number field and let $\phi(z) = z^2$. 
Fix a  non-archimedean place $v\in M_K$ such that $F_v$ 
is the unit disc in $\C_v$, and let $\ell$ be the residue characteristic of 
$\C_v$. Let $S_n = {\mathbb \mu}_{p^n}$ be the set of $p^n$-th
roots of unity. 
Clearly $\{S_n\}$ satisfies the conditions required in 
Theorem~\ref{thm2}. Hence, by Theorem~\ref{thm2} the sequence of sets
$\{ S_n \}$ is pseudo-equidistributed with respect to $(F_v, \hhat_{\phi,v})$. 
We can see this concretely as follows.

If $p \ne \ell$ then distinct elements are in distinct residue
classes, and it is easy to see for such a sequence of finite sets that
$\limsup_{n\to\infty} \d(S_n) = 1$, which is the transfinite diameter
of the unit disc. 

However, if $p = \ell$ then for all $\zeta \in 
{\mathbb \mu}_{p^n}$ we have $\zeta \equiv 1 \pmod{v}$. That is,
for all $n\ge 1$, $S_n$ is contained in the open unit disc centered 
at 1. A computation similar to that in
Remark~\ref{capacityofadisk}~(1) now shows that
$\limsup_{n\to\infty} \d(S_n) = 1$, which is not quite as
obvious as in the case $p \ne \ell$. 
\end{example}

\subsection{Equidistribution: archimedean case.}
% This section gives the proof of Proposition~\ref{balancedprop}

Our goal in this section is to prove Proposition~\ref{balancedprop}.
Before we do this, however, we need an easy lemma. 
Recall that a sequence of measures $\nu_n$ on a metric space $M$ is said to be 
{\it  tight} if for every $\epsilon > 0$, there exists a compact
set $K_{\epsilon} \subseteq M$ such that $\nu_n(M\setminus K_{\epsilon})
< \epsilon$ for all $n$ sufficiently large. Given $F \subset M$ and 
$\delta > 0$, we define 
\[
F_\delta := \{ z \in M \; : \; {\rm dist}(z,F) \leq \delta \}.
\] 

Suppose $ S =  \{ S_n \}$ is a sequence of finite subsets of
$M$. We say that $S$ is $F$-{\it tight} if given $\delta,\epsilon >
0$, there exists $N$ such that for
$n\geq N$, at least $(1-\epsilon)(\# S_n)$ of the elements of $S_n$ lie
in $F_\delta$.  

\begin{lemma}
\label{balancedlemma}
Let $\phi \in \C[z]$ be a polynomial with filled Julia set $F \subset
\C$.  Let $\lambda$ be a logarithmic distance function for $F$, and 
let $ S =  \{ S_n \}$ be a sequence of
finite subsets of $\C$ which is $\lambda$-taut.
Then $S$ is $F$-tight.
\end{lemma}

\begin{proof}
Let $\delta > 0$, and let $V_\delta$ be the complement of $F_\delta$.
Recall that the function $\lambda$ is continuous, nonnegative, and zero exactly on $F$, and that
$\lambda(z) \to \infty$ as $z \to \infty$.
By the local compactness of $\C$, it therefore follows that
$\lambda$ is bounded below by a positive
constant $C_\delta$ on $V_\delta$.
It follows immediately from this and the fact that
$\frac{1}{\# S_n} \sum_{z \in S_n} \lambda(z) \to 0$
that the sequence $\{ S_n \}$ is $F$-tight.
\end{proof}

We  now give the proof of Proposition~\ref{balancedprop}.

\begin{proof}
As usual, we let $N_n$ be the cardinality of $S_n$.  
Let $\mu_F$ denote the equilibrium measure for $F$, and let 
$\delta_n := \frac{1}{N_n} \sum_{z \in S_n} \delta_z$.
Since the sequence $\{ \delta_n \}$ is $F$-tight, it is in particular
tight in the usual sense, and so by Prohorov's theorem 
(see \cite{Parthasarathy})
it has a weakly convergent subsequence.
%\helpme{Ooops..., this sentence does not carry over to the 
%non-Archimedean case since $F_{\epsilon}$ is not compact even
%$F$ itself is compact.}
We may therefore assume that there is a probability
measure $\mu$ such that $\delta_n$ converges weakly to $\mu$, and it
suffices to show that $\mu = \mu_F$.  As a first step, we note that it
follows easily from the $F$-tightness of $\{ \delta_n \}$
that $\mu$ is supported on $F$.

Recall that the equilibrium measure on $F$
is characterized among probability measures
supported on $F$ as having
the smallest possible energy, which is precisely $-\log c(F)$.
% (The energy of a measure $\nu$ is by definition
% $-\int\int \log | x - y| d\nu(x) d\nu(y) $.)
Let $\Delta$ denote the diagonal in $\C \times \C$.
Then it suffices to prove the inequality

\begin{equation}
\label{eqn1}
\limsup_{n\to\infty}
\int\int_{\C \times \C - \Delta} \log |x - y| d\delta_n(x) d\delta_n(y)
\leq
\int\int_{\C \times \C} \log |x - y| d\mu(x) d\mu(y),
\end{equation}
because we know the left-hand side is equal to
\[
\log \left( \limsup_{n\to\infty} \prod_{x,y \in S_n, x\neq y}
|x - y|^{\frac{1}{N_n(N_n - 1)}} \right),
\]
which equals $\log c(F)$ since $S$ is a generalized Fekete sequence.

\smallskip

% To prove (\ref{eqn1}), 

Let $\epsilon >0$ be given, and set $S' := \{ S_n' \}$, 
where $S_n' := \{ z \in S_n \; : \; z \in F_\epsilon  \}$.
Then Lemma~\ref{balancedlemma} says that $S'$ is a full subsystem of $S$.
By Lemma~\ref{technicallemma}, it therefore suffices to show that

\begin{equation}
\label{eqn2}
\limsup_{n\to\infty}
\int\int_{F_\epsilon \times F_\epsilon - \Delta_\epsilon} \log |x - y| d\delta_n(x) d\delta_n(y)
\leq
\int\int_{\C \times \C} \log |x - y| d\mu(x) d\mu(y),
\end{equation}
where $\Delta_\epsilon$ is the diagonal in $F_\epsilon \times F_\epsilon$.

This inequality follows in an essentially formal manner from our assumptions and some
general measure theory.  Specifically, for a fixed $M$ we have
\begin{equation}
\label{eqnsmalldiag}
\int\int_{\Delta_\epsilon} \max \{ -M, \log |x - y| \} d\delta_n(x)
d\delta_n(y) = O(\frac{1}{N_n})
\end{equation}
by the definition of $\Delta_\epsilon$, and therefore
\[
\begin{array}{ll}
\limsup_{n\to\infty}
\int\!\int_{F_{\epsilon} \times F_{\epsilon} - \Delta_{\epsilon}
} \log |x - y|
d\delta_n(x) d\delta_n(y) &
\\
\le    \lim_{M\to \infty}
\limsup_{n\to\infty} \int\!\int_{F_{\epsilon} \times F_{\epsilon} -
  \Delta_{\epsilon}} 
\max\{-M, \log |x - y| \} d\delta_n(x) d\delta_n(y) &
\\
\indent\indent\indent\indent\indent\indent
\indent\indent\indent\indent\indent\indent
\indent\indent\indent\indent\indent
( {\rm since \;} (*) \le \max\{-M, (*) \} ) &
\\
=  \lim_{M\to \infty}
\limsup_{n\to\infty} \int\!\int_{F_{\epsilon} \times F_{\epsilon} }
\max\{-M, \log |x - y| \} d\delta_n(x) d\delta_n(y) & 
\\
\indent\indent\indent\indent\indent\indent
\indent\indent\indent\indent\indent\indent
\indent\indent\indent\indent\indent
( {\rm by \;} (\ref{eqnsmalldiag}) ) & 
\\
= \lim_{M\to \infty}  \int\!\int_{F_{\epsilon} \times F_{\epsilon}}
\max\{-M, \log |x - y| \} d\mu (x) d\mu(y) &
\\
\indent\indent\indent\indent\indent\indent
\indent\indent\indent\indent\indent\indent
\indent\indent\indent\indent\indent
(\delta_n \to \mu {\rm \; weakly}) & 
\\
= \int\!\int_{F_{\epsilon} \times F_{\epsilon}}
\log |x - y|d\mu (x) d\mu(y) &
\\
\indent\indent\indent\indent\indent\indent
\indent\indent\indent\indent\indent\indent
\indent\indent\indent\indent\indent
{({\rm monotone \; convergence \; theorem})} &
\\
=   \int\!\int_{\C\times\C} \log |x - y|d\mu (x) d\mu(y). & 
\\
\indent\indent\indent\indent\indent\indent
\indent\indent\indent\indent\indent\indent
\indent\indent\indent\indent\indent
{ ({\rm supp}(\mu) \subseteq F) } & 
\\
\end{array}
\]

%\helpme{Hmmm... why do we even need to replace S by S'?  I thought
%  there was a reason for it but at the moment I can't see what it is...}

% The inequality in the above chain comes from the obvious fact that
% $\log |x - y| \le \max\{-M, \log |x - y| \}$.
% The first equality is clear also as adding the diagonal $\Delta$
% to the integral will not affect the limit. The second equality
% follows from our assumption that $\delta_n$ converges to
% $\mu$ weakly, the third equality follows from the monotone
% convergence theorem, and the last equality comes from the fact
% that $\mu$ is supported on $F$.  

\end{proof}

%\helpme{Add a section here about $p$-adic equidistribution\ldots}
\subsection{Equidistribution: non-archimedean case.}
\label{padic-equidistribution}

Our goal is now to prove a non-archimedean analogue of
Proposition~\ref{balancedprop}. Let $L$ be an algebraically closed
non-archimedean local field. 
If we try to transport arguments used in
the archimedean case directly, we face some difficulties in the fact that the 
 non-archimedean field $L$  is not locally compact. 

For example, the analogue of Lemma~\ref{balancedlemma} is not in
general true over a non-archimedean field. 
This is clear if $F \subset L$ is the unit disc and 
$\lambda(z) = \log^+|z|$, since
$F_\delta = F$ for all $\delta \leq 1$. 
Another counterexample which may be illustrative is the following.

\begin{example}
\label{localheightexample}
Let $\phi(z) = (z^3 - z^2)/p$ be as in
Example~\ref{dynamicsexamples}(3).
Then if $|z|_p \leq \frac{1}{p}$, we have $\hhat_{\phi,p}(z) = 0$, and
if $|z|_p > 1$, we have $\hhat_{\phi,p}(z) = \log |z|_p + \frac{1}{2} \log p $.
For $\frac{1}{p} < |z|_p \leq 1$, there is no simple formula for the
canonical local height.  One can merely say that $\hhat_{\phi,p}(z)$ will
be zero if $z \in F$, and otherwise we will have $|\phi^n(z)|_p > 1$
for some $n$, in which case $\hhat_{\phi,p}(z) = 
\frac{1}{d^n} (\log |\phi^n(z)|_p + \frac{1}{2} \log p) $.

Recall from Example~\ref{dynamicsexamples}(3) that $F$ 
is contained in the union of the disc 
$ \{ v(z) \geq 1 \}$ and the sets $C_n := \{ z \in \C_p \; : \;
v(z) = 1 - \frac{1}{2^n}  \} $  
for $n = 0,1,2,\ldots$.  
Choose a rational number $0 < \epsilon_n < \frac{1}{2^n}$, and choose
$z_n \in \C_p $ with $v(z_n) = 1 - \frac{1}{2^n} - \epsilon_n$.  
% Then for all $z \in F$, we have $v(z - z_n) = \min \{ v(z),v(z_n)
% \} < 1$, so that $|z - z_n| > \frac{1}{p}$. 
If we set $y_n := \phi^n(z_n)$, then $v(y_n) = -2^n \epsilon_n$
and therefore 

\[\hhat_{\phi,p}(z_n) = \frac{1}{3^n} \hhat_{\phi,p}(y_n)
                      = \frac{1}{3^n} (2^n \epsilon_n \log p +
                      \frac{1}{2} \log p).
%  < \frac{\log p}{2\cdot 3^{n-1}}.
\]
Now it's clear that $\hhat_{\phi,p}(z_n) \to 0$ while the sequence
$\{z_n\}$ keeps a positive distance away from $F$ since
${\rm dist}(z_n,F) > \frac{1}{p}$.
\end{example}

%\begin{remark}
%\label{weirdcounterexample}
%The analogue of Lemma~\ref{balancedlemma} is not true if we replace
%$\C$ by $\C_p$ for some prime number $p$.  This is clear, for example, if $F
%\subset \C_p$ is the unit disk and $\lambda(z) = \log^+|z|$, since
%$F_\delta = F$ for all $\delta \leq 1$.  Another counterexample is
%furnished by the sequence $z_n$ from Example~\ref{localheightexample}:
%for all $z \in F$ we have  $v(z - z_n) = \min \{ v(z),v(z_n)
%\} < 1$, but we also have
%$ 0 < \hhat_{\phi,p}(z_n) < \frac{\log p}{2\cdot 3^{n-1}}$.  Therefore
%${\rm dist}(z_n,F) > \frac{1}{p}$ for all $n$ and $\hhat_{\phi,p}(z_n) \to 0$.

% When there was a part (ii) to the lemma, I had written:
% Similarly, one can construct a counterexample to (ii) by considering the
% polynomial $\psi(z) = p^3 z^3 - pz$, which is conjugate to $\phi$ via
% the map $z \mapsto p^2z$, and looking at the corresponding sequence $z_n' := p^{-2}z_n$.

%\begin{remark}
%On the other hand, it seems likely that the nonarchimedean analogue of
%Lemma~\ref{balancedlemma} is true under additional hypotheses.  Is it true, for
%example, if the filled Julia set of $\phi$ is compact?
%\end{remark}
%On the other hand, for non-Archimedean $L$ if the filled Julia set
%$F\subset F$ is compact then there is  a non-Archimedean analogue
%of Lemma~\ref{balancedlemma}(see Corollary~\ref{p-adicbalancedlemma}).
%To make our exposition uniform, we postpone proving this fact until
%Section~\ref{padic-equidistribution}.  
%\end{remark}

On the other hand, if the filled Julia set
$F\subset L$ is compact, then we will see that there is  a non-archimedean analogue
of Lemma~\ref{balancedlemma}.
We begin with the following result:

\begin{lemma}
\label{p-adiccompactJulia}
Let  $\phi(z) \in L[z]$ be a polynomial of degree $d \ge 2$. 
Assume that the filled Julia set $F$ of $\phi$ is compact.
Then for each $\delta \in |L^{\ast}|$, there exists
a positive integer $N$ such that 
$F_n \subset F_{\delta}$ for all $n \ge N.$  
\end{lemma}

\begin{proof}
Our assumption implies that the filled Julia set equals the Julia set 
of $\phi$,  for otherwise there would be a closed disc totally contained in 
$F$ by \cite[Theorem~5.1 (3),(4)]{Benedetto-Component}, and hence
$F$ would not be compact.  

Let $\delta \in |L^{\ast}|$ be given. For $\omega \in L$,
let $D_{\delta}(\omega) = \{ z \in L : |z - \omega| < \delta \}$
denote the open disc centered at $\omega$ of radius $\delta$. Then
the collection $\{D_{\delta}(\omega) \; : \; \omega \in F \}$ of open
discs gives an open covering of $F$. Since $F$ is compact, there 
exist $\omega_1, \ldots, \omega_{m(\delta)} \in F$ such that 
$F \subseteq U_{\delta} = \cup_{i=1}^{m(\delta)} D_{\delta}(\omega_i)$. 

Since $F$ is the Julia set of $\phi$, for each $i$ there exists
a positive integer $N_i$ such 
$\phi^{N_i}(D_{\delta}(\omega_i)) \supseteq F_0$. Take 
$N = \max \{N_1, \ldots, N_{m(\delta)}\}$. Then 
$\phi^N(D_{\delta}(\omega_i)) \supseteq F_0$ for all $i$. 
We claim that $F_N \subseteq  U_{\delta}$.

Assume to the contrary that $F_N$ is not contained in $U_{\delta}$.
We note that $F_N$ itself is a disjoint union of closed discs 
$F_N = \cup_{j} B_{N,j}$ where $\phi^N (B_{N,j}) = F_0$. As the valuation
on $L$ is ultrametric, any two discs are either disjoint or one
contains the other.
Therefore, there must be some $B_{N, l}$ which is disjoint
from $U_{\delta}$. This is impossible as $B_{N,l}$ contains the
inverse image of some point of $F$ under $\phi^N$ and hence has non-empty intersection
with $F$. 

Our claim  now shows that $F_N \subseteq F_{\delta}$.  
On the other hand, we have $F_N \supseteq F_{N+1} \supseteq \cdots$. 
Now the assertion of our proposition is clear.    
\end{proof}

The following
analogue of Lemma~\ref{balancedlemma}
suffices in order to prove our non-archimedean equidistribution theorem.

\begin{prop}
\label{p-adicbalancedlemma}
Let $\phi \in L[z]$ be a polynomial with compact filled Julia set 
$F \subset L$.  
Let $\lambda_{\phi}$ be the canonical local height associated
to $\phi$ as defined in Section~\ref{localheightsection}, and 
let $ S =  \{ S_n \}$ be a sequence of
finite subsets of $L$ which is $\lambda_{\phi}$-taut.
Then $S$ is $F$-tight.
\end{prop}

\begin{proof}
As in the proof of Lemma~\ref{balancedlemma}, it suffices 
to show that $\lambda_{\phi}$ is bounded below by a positive
constant outside $F_{\delta}$ for any $\delta \in |L^{\ast}|$. 

Recall that if $F_0 = \{ |z| \leq R \}$ with $R$ sufficiently large, then
$|\phi^m(z)| = |a(m) z^{d^m}| > \max(1, R)$ for all $z \not\in F_0$ where 
$a(m) := a_d^{1 + d + d^2 + \cdots + d^{m-1}}$.  
It follows that $\lambda_{\phi}(z) = \log |z| - \log c(\phi)$ 
for all $z\not \in F_0$, where $c(\phi) = |a_d|^{\frac{-1}{d-1}}$ as
usual.  From this formula, it is clear that
$\lambda_{\phi}$ is bounded below by a positive constant outside
$F_0$. 

Let $\delta \in |L^{\ast}|$ be given. By 
Lemma~\ref{p-adiccompactJulia}
we know that there is a positive integer $N$ such that 
$\phi^n(z) \not \in F_0$ for $z \not\in F_{\delta}$ and $n \geq N$. 
As 
\[
\lambda_{\phi}(z) =  (1/d^N) \lambda_{\phi}(\phi^N(z)),
\]
for all $z \not\in F_{\delta}$, our assertion now follows in general.
\end{proof}

% \begin{remark}
% The arguments for the two assertions above works for general 
% ultrametric valued field.
% \end{remark}

Recall from \cite[Theorem~4.1.22]{Rumely} that for a compact 
subset $E \subseteq L$ of positive transfinite diameter, there exists a unique equilibrium 
probability measure $\mu_E$ supported on $E$.  
We now prove Proposition~\ref{nonarch-balancedprop}.

\begin{proof}(of Proposition~\ref{nonarch-balancedprop}):

As in the proof of Proposition~\ref{balancedprop}, 
we let $\delta_n := \frac{1}{N_n} \sum_{z \in S_n} \delta_z$.

In the first step, we need to show that the sequence
of measures $\{\delta_n\}$ has a weakly convergent 
subsequence. Equivalently,  we will  show
that for any given $\epsilon, \delta > 0$ there exists
a set $F_{\epsilon, \delta}$ which is the union of a finite
number of discs of radius $\delta$ such that 
$\delta_n(F_{\epsilon, \delta}) > 1 - \epsilon$ for 
all $n$ (see \cite[p. 49]{Parthasarathy}). 

By  Corollary~\ref{p-adicbalancedlemma}  
the sequence $S = \{S_n\}$ is $F$-tight, for any given
$\epsilon, \delta$ it is immediate from the definition
that there exists a integer $N$ such that 
$\delta_n (F_\delta) > 1 - \epsilon$ for all   
$n \ge N$. We note that $F_{\delta}$ is already
a finite union of discs of radius $\delta$. 
Let $F_{\epsilon,\delta}$ be 
the union of $F_{\delta}$ and finitely
many discs centered at points of the set $\cup_{i < N}\,S_i$.
Then, clearly we have $\delta_n(F_{\epsilon,\delta}) > 1 - \epsilon$
for all $n$.
Thus, we may  assume that there is a probability
measure $\mu$ such that $\delta_n$ converges weakly 
to $\mu$ and show that $\mu = \mu_{F}$. 

The remaining arguments
are the same as  in the proof of Proposition~\ref{balancedlemma}
and we will not repeat  them  here.
 \end{proof}

\section{Applications}
\label{Applications}

\subsection{Lower Bounds for Canonical Heights}

% In this paragraph, we will restrict our attention to the case where
% $K$ is either a number field or a function field of transcendence 
% degree one over a finite field and 
% present a version of non-archimedean equidistribution
% theory. Let $v\in M_K$ be a non-archimedean place of $K$. 

Let $v\in M_K$ be a non-archimedean place. 
We say an algebraic extension $E$ of $K$ is totally
$v$-adic of type $(e, f)$ (over $K$)  if for any 
embedding $E\hookrightarrow \C_v$, the image of $E$ is 
contained in a finite extension of $K_v$ whose ramification degree
and residue degree are bounded by $e$ and $f$ respectively.
Likewise, an $\alpha \in \Kbar$ is said to be totally 
$v$-adic of type $(e,f)$ if
$K(\alpha)$ is totally $v$-adic of type $(e,f)$.

% The following is wrong, because the composite doesn't preserve e and
% f:
%
% We say that $E$ is maximally $v$-adic extension of type $(e,f)$ over $K$ 
% if  $E$ is the composite of all
% the $v$-adic extensions of type $(e,f)$ over  $K$. 
% Note that by our definition, if $E$ is a $v$-adic extension of
% type $(e,f)$ over $K$ then its normal closure is also a 
% $v$-adic extension of type $(e,f)$ over $K$.

\begin{theorem}
\label{padicSchinzel}
(c.f. \cite[Theorem 2.]{Bombieri-Zannier})  Let 
$\phi(z) = a_0 + a_1 z + \cdots + a_d z^d \in K[z]$ and let $v\in M_K$ 
be  a place of good reduction for $\phi$.
%non-archimedean place of $K$ such that $v(a_i) \ge 0 $ for all
%$i$ and $v(a_d) = 0$.
%\helpme{Maybe define good reduction and say ``a place of good
%  reduction for $\phi$'' instead.}
Let $e, f$ be positive integers. 
Then there exists a 
positive constant $C = C(e, f, v)$ such that $\hhat_{\phi}(\alpha) \ge C$
for all but finitely many $\alpha$ which are totally $v$-adic of type $(e,f)$. 
In particular, 
there are only finitely many preperiodic points of $\phi$ 
which are totally $v$-adic of type $(e,f)$.   
\end{theorem}   

\begin{remark}
If $K$ is a number field and $\phi(z) = z^d $ for some $d \ge 2$, then $h_{\phi}(\alpha) =
h(\alpha)$ is just the absolute logarithmic height of
$\alpha$.  In this special case, Theorem~\ref{padicSchinzel} can
be viewed as a qualitative version of \cite[Theorem~2]{Bombieri-Zannier}.
\end{remark}

\begin{proof}
Let $E$ be the 
compositum of all finite extensions of $K_v$ whose ramification indices and
residue degrees are bounded by $(e,f)$. 
Thus, $E$ contains all $\alpha \in \Kbar$ which are totally $v$-adic 
of type $(e,f)$. Note that  $E$ is a finite extension over 
$K_v$, since there are only finitely many extensions of $K_v$ with prescribed 
ramification indices and residue degrees. Let $e', f'$ be the 
ramification index and the residue degree, respectively, of $E$ over $K_v$. 
 
It suffices 
to show that there does not exist a sequence $ \{ S_n \}$ in $E$ of distinct
finite {\it Galois-stable} subsets of $\Kbar$
such that $\lim_{n\to \infty} \hhat_\phi (S_n) = 0$.
Assume to the contrary that there is such a sequence 
$ S = \{ S_n \}$. Then $S$ is $\hhat_{\phi,v}$-taut. 
Let $S_n' = \{\alpha \in S_n \; : \; |\alpha|_v \le 1\}$. 
We claim that $S_n'$ is full subsystem  of $S_n$, i.e., that 
$\lim_{n\to \infty} \# S_n' / \# S_n = 1$. 
Equivalently, we need 
to show that $\lim_{n\to \infty} \# T_n / \# S_n = 0$, where 
$T_n = S_n \setminus S_n'$. 

By hypothesis, we have 
$\hhat_{\phi,v} (z) = \log^+ |z|_v$, so that
$$
\hhat_{\phi,v}(S_n) = \frac{1}{\# S_n} \sum_{\alpha\in T_n} \,
\log^+ |\alpha |_v .
$$

Let $p$ be the residue characteristic of $\C_v$.
As the ramification degree of $E$ over $K_v$ is $e'$,
we see that $\log^+ |\alpha|_v \ge (1/e') \log p$ for 
$\alpha \in T_n$. It follows that
$\hhat_{\phi,v}(S_n) \ge \frac{\# T_n}{e'\, \# S_n}  \log p$. If 
$\# T_n / \# S_n \not\to 0$ then $S_n$ is not $\hhat_{\phi,v}$-taut
which contradicts our assumption. Hence, $S_n'$ is a full
subsystem of $S_n$ as desired. 

By Lemma~\ref{technicallemma}, $\limsup \d(S_n) \leq \limsup \d(S_n').$
However, we have $\limsup \d(S_n') \le c(\Ocal_E)$, where $\Ocal_E$ is
the ring of integers of $E$. 
We have $c(\Ocal_E) = p^{-1/e'(p^{f'} - 1)} < 1$
as computed in Example~\ref{computetrasfinitediameters}(2). This shows that
$S = \{S_n\}$ is not pseudo-equidistributed and contradicts 
Theorem~\ref{thm2}.   
\end{proof}

% \helpme{Maybe we should compare the above result  with
% Bombieri-Zannier's. Or, improve the above result as close 
% as possible to theirs.} 

Let $p$ be a rational prime. 
Following \cite{Bombieri-Zannier}, we call an
algebraic number $\alpha$ {\em totally $p$-adic} if the prime $p$ splits
completely in $\Q(\alpha)$. This is equivalent, in the above
terminology, to saying that $\alpha$
is totally $p$-adic of type $(1,1)$ over $\Q$. As a special case of 
Theorem~\ref{padicSchinzel}, we have the following $p$-adic analogue
of Schinzel's theorem \cite{Schinzel} on heights of totally real
algebraic numbers:

\begin{cor}
(c.f. \cite[Example~1.]{Bombieri-Zannier})  
Let $\phi(z) = a_0 + a_1 z + \cdots + a_d z^d \in \Q[z]$, and let $p$
be a prime where $\phi$ has good reduction.  Then there exists a 
constant $C>0$, depending only on $p$ and $\phi$, such that 
$\hhat_{\phi}(\alpha) \ge C$ for all but finitely many 
totally $p$-adic algebraic numbers $\alpha$. 
\end{cor}

We can also give an archimedean generalization of Schinzel's theorem.
% , though our generalization is weaker in the
% sense that the lower bound of the canonical height is not explicit. 
We use the fact (see \cite[Proof of Theorem~6.5.8]{Ransford}) that if $F$ is the filled Julia set of
a polynomial $\phi \in \C[z]$ of degree at least 2, then 
the support of $\mu_F$ is precisely the Julia set of $\phi$.

\begin{prop}
\label{Schinzelresult}
Let $K$ be a number field and let $\phi \in K[z]$ of degree $d \ge 2$. 
Assume that  the complex Julia set $J_\phi$ of $\phi$ is not contained
in the real line. Then there exists a constant 
$C>0$, depending only on $\phi$, such that for all but finitely
many totally real algebraic numbers $\alpha$, 
we have $\hhat_{\phi}(\alpha) \geq  C$. 
\end{prop} 

\begin{proof}
Assume to the contrary that there exists 
sequence of totally real algebraic numbers $\{\alpha_n\}$ such 
that $\hhat_{\phi}(\alpha_n) \to 0$.  By Corollary~\ref{thm2cor},
$\{\alpha_n\}$ is equidistributed with respect to $\mu_F$, where $F$ is
the filled Julia set of $\phi$. 
This implies that $\mu_F$ is supported on a subset of the real 
line.  On the other hand, the support of $\mu_F$ is equal to the Julia set
of $\phi$, which by assumption is not totally contained in the real line.
We thus have a contradiction.
\end{proof}

% \helpme{For non-Archimedean place, assume that  the Julia set is not
% compact. Let $F(L)$ denote the intersection of the filled Julia set
% with a discrete valued subfield $L$. 
% Is it true that the transfinite diameter of $F(L)$ is strictly
% less than the transfinite diameter of the full filled Julia set?} 
%
% I don't know: probably yes.

For the rest of this section, we assume $K$ is a number field, and we
fix an embedding of $K$ into the complex numbers.
For the next application of our main result, we are concerned with 
an analogue of a result
of Zhang \cite{Zhang, Zagier} for the usual height on $\Qbar$ 
in the setting of dynamical canonical heights.
We want to know how small the canonical height
of algebraic points on a line in ${\mathbf A}^2$ can be.  

\begin{lemma}
\label{affinetransform}
Let $ \{ S_n \}$ be a sequence of distinct,
finite, Galois-stable subsets of $\Kbar$
such that $\lim_{n\to \infty} \hhat_\phi (S_n) = 0$. 
Let $m(z) = \alpha(z + \beta)$ be an affine 
linear transformation defined over $K$, with $\alpha \neq 0$, and
suppose $\lim_{n\to \infty} \hhat_\phi (m(S_n)) = 0$.
Then \\
(1) $\alpha$ is a root of unity, and 
\\
(2) $m (J_{\phi}) = J_{\phi}$ where $J_{\phi}$ denotes
the complex Julia set of $\phi$.  
\end{lemma}

\begin{proof}
It follows from Theorem~\ref{thm2} that both sequences
$\{ S_n \}$ and $\{ m(S_n) \}$ are 
pseudo-equidistributed with respect to the pair $(F_v,\hhat_{\phi,v})$
for each $v\in M_K$. Therefore, for any  $v \in M_K$
\[
\lim_{n\to \infty} \d_v(S_n) = c_v(\phi) = \lim_{n\to \infty} 
\d_v (m (S_n)).
\]

Clearly, $\d_v(m (S_n)) = |\alpha|_v \d_v (S_n)$.
Thus, since $c_v(\phi) > 0$, we have $|\alpha|_v = 1$ for all $v \in
M_K$, and it follows that $\alpha$ is a root of unity.  This
proves (1).

For the proof of (2), note that by Corollary~\ref{thm2cor}
both sequences $\{ S_n \}$ and $\{ m(S_n) \}$ are 
equidistributed with respect to $\mu_F$, where $\mu_F$ is
the equilibrium measure on the filled Julia set $F$. 
It follows easily from this that $\mu_F = \mu_{m(F)}$.  But the
support of $\mu_F$ is $J_\phi$ and the support of $\mu_{m(F)}$ is
$m(J_\phi)$, so $J_\phi = m(J_\phi)$ as desired.

% Previous argument:
% \[
% \delta_{S_n} \Rightarrow \mu_F \quad{\text and} \quad 
% \delta_{m(S_n)} \Rightarrow \mu_F
% \]
% Moreover, we also have $\delta_{m(S_n)} \Rightarrow \mu_{m(F)}$. 
% 
% On the other hand, as $\mu_F$ is 
% supported on the Julia set $J_{\phi}$, 
% if $m(J_{\phi}) \ne J_{\phi}$ then we can choose  any
% point $x \in m(J_{\phi}) \setminus J_{\phi}$ and consider an 
% open disk $D$ containing $x$. but disjoint from
% $J_{\phi}$.  Note that we must have 
% $x$  in the Fatou set of $\phi$. Therefore, we 
% can take  $D$ to be  disjoint from
% $J_{\phi}$. 
%
% Then, $\int_{D}\, \delta_{m(S_n)} \to \int_{D} d\mu_F = 0$ 
% while $\int_{D}\, \delta_{m(S_n)} \to \int_{D} d\mu_{m(F)} > 0$
% which is a contradiction. Therefore, we must have $m(J_{\phi}) 
% = J_{\phi}$. 
\end{proof}

\begin{remark}
\label{translation}
If $\alpha = 1$ then we must have $\beta = 0$ in
Lemma~\ref{affinetransform}.
This follows immediately from that fact that $J_\phi$ is a bounded
subset of $\C$, since if $\beta \ne 0$ then clearly $J_{\phi} + \beta 
\ne J_{\phi}$, which contradicts (2) above.   
\end{remark}

We see that under the hypotheses of Lemma~\ref{affinetransform}(2), 
% if $\{S_n\}$ and $\{m(S_n)\}$ 
% are both sequences of distinct finite Galois-stable subsets of $\Kbar$ 
the affine transformation $m$ must be in the symmetry group of the Julia set. This 
group has been studied in \cite{Beardonsymmetry} (see also the
references cited in that paper). 
In order to state what we need from \cite{Beardonsymmetry}, we 
introduce some notation.  We define the {\it centroid} of 
$\phi(z) = a_d z^d + a_{d - 1} z^{d-1} + \cdots + a_1 z + a_0 
\in \C[z]$ to be the complex number $\zeta = - a_{d-1}/(d\, a_d)$.
We also let $\Sigma(\phi)$ denote the set of all affine linear transformations
$\sigma$ on the complex plane such that $\sigma(J_{\phi}) = J_{\phi}$,
and call $\Sigma(\phi)$ the set of symmetries of $J_\phi$.

\begin{theorem}
\label{Juliasymmetry}
\begin{itemize}
\item[1.] The set $\Sigma (\phi)$ of symmetries of $J_{\phi}$ consists of the
rotations about the centroid of $\phi$.
\item[2.] If $\sigma$ is any affine transformation, then 
$\sigma \in \Sigma (\phi)$ if and only if $\phi\circ \sigma = \sigma^{d}
\circ \phi$. 
% I deleted the rest because I don't think it's as relevant for us as
% the first part of the theorem.
%
% and if $a_{d-1} = 0$ 
% and $\Sigma (\phi)$ is finite, then the order of $\Sigma (\phi)$ 
% is the largest integer $n$ such that $\phi$ can be written in the 
% form $\phi(z) = z^r \psi(z^n)$ for some polynomial $\psi$. 
\end{itemize}
\end{theorem}

\begin{proof}
See Theorem 5 and Lemma 7 of \cite{Beardonsymmetry}.
\end{proof}

Note that part 2 of Theorem~\ref{Juliasymmetry} implies 
that a point $Q$ is preperiodic for $\phi$ 
if and  only if $\sigma(Q)$ is.  
% It will be interesting to see if the 
% converse of Lemma~\ref{affinetransform}(2) is true. 
The only
case where $\Sigma(\phi)$ is infinite is when $\phi$ is conjugate
to a polynomial map of the form $z\mapsto z^n$ 
for some integer $n \ge 2$ \cite[Lemma 4]{Beardonsymmetry}.
%\helpme{What about polynomials conjugate to $z^n$?}

\medskip

We now consider a line $L$ defined by an equation of the form
$a x + b y = c \; (a, b, c \in K)$,  and we study the canonical height
of algebraic points $(x, y) \in L(\Kbar)$. Here, we embed 
${\mathbf A}^2$ into $\P^1 \times \P^1$ via
$(x, y) \mapsto ([x: 1], [y:1])$, we consider
$\Phi = (\phi(x), \phi(y))$ as a self map on $\P^1 \times \P^1$, and
we define the canonical height of an algebraic point $Q = (x,y)$ 
to be $\hhat_{\Phi}( (x, y) ) = \hhat_{\phi}(x) + \hhat_{\phi}(y)$. 
%As remarked above, if $\Sigma (\phi)$ is infinite then 
%$\phi(z) = z^n$ for some integer $n\ge 2$. In this case, 
%a lower bound for the canonical heights of algebraic 
%points on $L$ is given by Zhang's Theorem \cite{Zhang}
%(see also \cite{Zagier}).  
% We therefore need only consider the
% case where $\Sigma (\phi)$ is finite.

%Moreover, if 
If one of $a, b$ is zero than the closure of $L$ in
$\P^1 \times \P^1$  is of the form
$\{\alpha\}\times \P^1$ or $\P^1 \times \{\beta\}$,
and it is easy to see that $L$ contains sequence $\{S_n\}$ of Galois 
stable algebraic points with $\hhat_{\Phi}(S_n) \to 0$
if and only if $\alpha$ (respectively, $\beta$) is a 
preperiodic point of $\phi$.

\begin{prop}
\label{Bogomolov}
Assume that $\Sigma(\phi)$ is finite.
% a finite group of exponent $e$.
Let $a, b, c \in K$ with $a b \ne 0$, and assume either that
$ - a/b$ is not a root of unity 
% (More specifically, an $e$-th root of unity)
or that $c/a \ne \zeta $, where $\zeta$ is the centroid of $\phi$. 
Let  $L\subset {\mathbf A}^2$ be the line
defined by the  equation $a x + b y =  c$. Then there exists 
a positive constant $C = C( a, b, c)$ such that 
the set $\{(x, y) \in L(\Kbar) \; : \; \hhat_{\Phi}((x,y)) :=
\hhat_{\phi}(x) + \hhat_{\phi}(y) < C \}$ is finite.
In particular, $L$ contains only finitely many points 
$(s, t)$ such that both $s$ and $t$ are preperiodic points of $\phi$.
\end{prop}

\begin{proof}
Assume to the contrary that there exist a sequence of points 
$\{(x_n, y_n) \in L(\Kbar)\}$ such that 
$\lim_{n\to \infty}\, \hhat_{\Phi}((x_n, y_n)) = 0$. Then we 
have two sequences $\{x_n\}$ and $\{y_n\}$ such that 
$\hhat_{\phi}(x_n) \to 0 $ and $\hhat_{\phi}(y_n) \to 0$. 

As $y_n = (\frac{-a}{b}(x_n - \frac{c}{a}))$, we see by combining 
Lemma~\ref{affinetransform} and Theorem~\ref{Juliasymmetry} 
that it is impossible for $\hhat_{\phi}(x_n) \to 0 $ and 
$\hhat_{\phi}(y_n) \to 0$ to hold simultaneously. 
\end{proof}

\begin{remark}
One can replace the assumption that $-a/b$ is not a root of unity in
the statement of Proposition~\ref{Bogomolov} with the stronger
hypothesis that $-a/b$ is not an $e$th root of unity, where $e$ is the
exponent of $\Sigma(\phi)$.
Moreover, whenever $\Sigma (\phi)$ is finite, one can show that 
the order of $\Sigma (\phi)$ 
is the largest integer $n$ such that $\phi$ can be written in the 
form $\phi(z) = z^r \psi(z^n)$ for some polynomial $\psi$ (see \cite{Beardonsymmetry}).
\end{remark}

This result motivates the following theorem, which is closely related
to Conjecture 2.5 of \cite{Zhang2}.

\begin{theorem}
\label{Bogomolov-Zhang}
% Suppose that $\Sigma(\phi)$ is finite, and
Let $L/K$ be any line in ${\mathbf A}^2$.  Then the following are
equivalent:
\begin{itemize}
\item[(a)] $L$ is a preperiodic line under $\Phi$.
\item[(b)] $L$ is defined by an equation of the form $x = \alpha$ 
or $y = \alpha$, with $\alpha$ a preperiodic point of $\phi$, or 
$y = \sigma(x)$ with $\sigma \in \Sigma(\phi)$.
\item[(c)] $L$ contains infinitely many preperiodic points of $\Phi$.
\item[(d)] $L$ contains an infinite sequence of algebraic points
  $(x_n,y_n)$ such that \\ $\hhat_\Phi(x_n,y_n) \to 0$.
\end{itemize}
\end{theorem}

\begin{proof}

We may assume without loss of generality that $\Sigma(\phi)$ is
finite, since otherwise $\phi(z)$ is conjugate to $z^n$ for some $n
\ge 2$ and the result follows from \cite{Zhang}. 

(a) implies (d): 
Since $L$ is preperiodic, there exist positive integers $k$ and $k'$
and a curve $C$ such that $\Phi^{k'} L = C$ and $\Phi^{k} C = C$.
The maps $\Phi^{k'} : L \to C$ and $\Phi^k : C \to C$ of affine algebraic 
curves extend to nonconstant
maps of projective algebraic curves, and therefore must be surjective.
It follows that for any point $(x_0,y_0) \in C$ 
and any $n \geq 1$, we can find 
a point $(x_n,y_n) \in L$ such that $\phi^{nk+k'}(x_n,y_n) = (x_0,y_0)$.  
We then have
$\hhat_\Phi(x_n,y_n) = d^{-nk-k'} \hhat(x_0,y_0) \to 0$ as $n \to
\infty$.
Moreover, 
if we choose $(x_0,y_0)$ so that it is not a 
preperiodic point for $\Phi^k$, then
the points $(x_n,y_n)$ are distinct.

(d) implies (b): If the equation of $L$ is not of the form  $x = \alpha$ or $y =
\alpha$, with $\alpha$ a preperiodic point of $\phi$, then it follows by
the same proof as that given for Proposition~\ref{Bogomolov} that 
$L$ is given by an equation $y = \sigma(x)$ with $\sigma \in \Sigma(\phi)$.

(b) implies (a): If $L$ is defined by $y = \sigma(x)$ with 
$\sigma \in \Sigma(\phi)$,
then $\Phi(L)$ is defined by $y = \sigma^d(x)$.  As $\sigma$ has
finite order, it follows that $L$ is preperiodic.  

(b) implies (c): Let $x_0$ be any preperiodic
point for $\phi$.  If $L$ is defined by $y = \sigma(x)$ with 
$\sigma \in \Sigma(\phi)$,
then since $\sigma$ preserves the set of $\phi$-preperiodic points,
it follows that $L$ must contain the $\Phi$-preperiodic point $(x_0,
\sigma(x_0))$.

(c) implies (d): Obvious.
\end{proof}

% \begin{remark}
% When $\phi(z) = z^n$ for some $n \ge 2$, then $\Sigma(\phi)$ is an
% infinite group. In this case, Theorem~\ref{Bogomolov-Zhang} is a 
% special case of Zhang's Theorem \cite{Zhang}. 
% \end{remark}

% \begin{remark}
% (1).
% If the line $L$ is defined by the equation $y = \sigma(x) = \lambda x + 
% \gamma$ where $\sigma (z) \in \Sigma (\phi)$ 
% then as remarked above 
% we have that $\sigma$ preserves the set of preperiodic points of 
% $\phi$. Thus, $L$ contains the set of points $\{ (x, \sigma(x))\mid
% \, \hhat_{\phi}(x) = 0\}$ which are points with zero canonical 
% heights. It is not clear if $\{ S_n \}$ with $\hhat_{\phi}(S_n) \to 0$
% then  $\{ \sigma (S_n)\}$  also has $\hhat_{\phi} (\sigma (S_n)) \to 0$. 
% \smallskip
% \\
% (2). If $L$ is given by $y = \sigma(x)$ for some $\sigma\in \Sigma(\phi)$
% as in (1), then we have $\Phi(L)$ is also a line defined by
% the equation $y = \sigma^d (x)$ by \cite[Lemma 7.]{Beardonsymmetry}
% already quoted. As $\sigma$ is of finite order, it follows that 
% $L$ is a {\em preperiodic line} under iterates of $\Phi$. It is interesting
% to know whether or not the converse is true. Namely, if $L$ contains
% infinitely many point $(x, y)$ with both $x, y$ are preperiodic
% points of $\phi$ then is it necessary that $L$ is preperiodic under 
% the iterates of $\Phi$. 
% 
% \end{remark}

\subsection{Generalized height functions and the Mandelbrot set}

\label{Mandelbrot}

In this section, 
% we formulate a more general notion of height, 
we prove a more general version of Theorem~\ref{thm2}.  As an application, we
attach a height function to the (adelic) Mandelbrot set and obtain a
corresponding equidistribution theorem.

First, we define generalized Green's functions.  Let $L$ be an
algebraically closed local field, let $F$ be a closed and bounded
subset of $L$, and let $\lambda$ be a logarithmic distance function for $F$.  
If $F$ is a lemniscate of the form $\{ z\in L \; : \; |P(z)| \leq R
\}$ for some polynomial $P \in L[z]$ of degree $d$ and some $R$ in the
value group of $L$, we call the function 
\[
g_F(z) = \frac{1}{d} \log^+ |\frac{P(z)}{R}|
\]
the {\it Green's function} attached to $F$.  
It follows from \cite[Proposition~4.4.1]{Rumelybook} that $g_F(z)$
is well-defined, independent of the choice of a particular defining
pair $(P(z),R)$.

%In the following, we'll call $F$ a {\it generalized lemniscate} set if
%there exists a descending sequence $F_1\supseteq F_2 \supseteq
%\cdots$ of lemniscates of the form $F_n = \{ z\in L \; \mid \;
%|P_n(z)| \leq R_n \}$ for some polynomial $P_n \in L[z]$ of degree
%$d_n$ and some positive number $R_n$ in the value group of $L$ such
%that  $\cap_{n=1}^\infty F_n = F$. 
%We say that $\lambda$ is a {\it generalized Green's function} for $F$
%if $F$ is generalized lemniscate and 
%\begin{itemize}
%\item $\lim_{n \to\infty} c(F_n) = c(F)$ (this does not automatically
%  follow from the previous condition -- see Remark~\ref{nesteddiskremark})
%\item $\lambda(z) = \lim_{n \to\infty} g_n(z)$, where $g_n = g_{F_n}$ is the
%  Green's function attached to $F_n$.
%\end{itemize}
More generally, we say that $\lambda$ is a {\it generalized Green's function} for $F$
if there exists a descending sequence $F_1\supseteq F_2 \supseteq
\cdots$ of lemniscates of the form $F_n = \{ z\in L \; : \;
|P_n(z)| \leq R_n \}$ for some polynomial $P_n \in L[z]$ of degree
$d_n$ and some positive number $R_n$ in the value group of $L$ such
that:
\begin{itemize}
\item $\cap_{n=1}^\infty F_n = F$
\item $\lim_{n \to\infty} c(F_n) = c(F)$ (this does not automatically
  follow from the previous condition -- see Remark~\ref{nesteddiskremark})
\item $\lambda(z) = \lim_{n \to\infty} g_n(z)$, where $g_n = g_{F_n}$ is the
  Green's function attached to $F_n$.
\end{itemize}

\begin{remark}
\label{Rumelyremark}
It follows from the results of \cite[Section 4.4]{Rumelybook}
(especially Theorem 4.4.4 and Lemma 4.4.7)
that if $F$ is {\it algebraically
  capacitable}, then there exists a unique generalized Green's
  function attached to $F$ called (naturally) the Green's
  function of $F$ (with respect to the point $\infty$).  As shown in \cite{Rumelybook}, examples of
  algebraically capacitable sets include all finite unions of compact sets and
  lemniscates in $L$.  However, we do not know, for example, whether
  or not the $v$-adic filled Julia set of a 
  polynomial $\phi$ is always algebraically capacitable.  This is why we have
  chosen to define generalized Green's functions.
\end{remark}

The following important fact is proved in \cite[Proposition 4.4.1(C)]{Rumelybook}:

\begin{lemma}
\label{greensfunctionlemma}
If $F_1$ and $F_2$ are lemniscates with $F_1 \subseteq F_2$, then
\[
g_{F_1}(z) \geq g_{F_2}(z)
\]
for all $z \in L$.
\end{lemma}

% \helpme{I think that if we don't use this lemma, then we can
%   still get by if we assume that the $g_n$'s converge {\it uniformly}
%   to $\lambda$.  Is that right?  In any case, I think it's better to
%   use the lemma from Rumely's book.}

We then have the following generalization of Lemma~\ref{keylemma}:

\begin{lemma}
\label{keylemmavariant}
Let $F$ be a closed and bounded subset of $L$, and suppose that
$\lambda$ is a generalized Green's function for $F$.
Suppose that $S_n$ is a sequence of finite subsets of $L$ with the property that
$N_n := \# S_n \to \infty$ and  $\lambda(S_n) \to 0$.
Then
\[
\limsup_{n\to\infty} \d(S_n) \leq c(F).
\]
\end{lemma}

\begin{proof}
Fix $\beta > 0$ in the value group of $L$, and set
$S_n' := \{ z \in S_n \; : \; \lambda(z) \leq \beta \}$.
Then as in the proof of Lemma~\ref{keylemma},
we see easily that $S' := \{ S_n' \}$ is a full subsystem of $S$.

Moreover, it follows from Lemma~\ref{greensfunctionlemma} and the fact
that $g_n \to \lambda$ that for all $m\geq 1$ we have
\[
\{ z \in L \; : \; \lambda(z) \leq \beta \} \subseteq 
F_{m,\beta} := \{ z \in L \; : \; g_m(z) \leq \beta \}.
\]

Therefore $S_n' \subset F_{m,\beta}$ for all $m$, so that 
\[
\limsup_{n\to\infty} \d(S_n') \leq c(\cap_{m=1}^\infty F_{m,\beta})
= c(F)e^\beta
\]
by Remark~\ref{PLremark} and the assumption that $c(F_n) \to c(F)$.

Letting $\beta \to 0$ now gives the desired result.
\end{proof}

Now let $K$ be a global field, and let $\F = \{ F_v \}$ be an adelic
set with respect to $K$ (see Section~\ref{Fekete-Chebyshev}).  
We assume that $\F$ also comes equipped with
an {\it adelic generalized Green's function}, by which we mean
a collection $\g := \{ \lambda_v \}$ consisting of a generalized Green's
function $\lambda_v$ for $F_v$ for each $v \in M_K$.

The height function $h_{\F} : \Kbar \to \R$ attached to $\F$ (or more
properly to the pair $( \F,\g )$) is then defined by
\[
h_{\F}(\alpha) := \frac{1}{\# S} \sum_{v \in M_K} \sum_{\alpha_i \in S}
n_v \lambda_v(\alpha_i),
\]
where $S$ is the smallest Galois-stable subset of $\Kbar$ containing $\alpha$.

In this context, we have:

\begin{theorem}
\label{thm2variant}
Let $K,\F,\g$ be as above, and assume that $c(\F) = 1$.
Let $ \{ S_n \}$ be a sequence of distinct
finite Galois-stable subsets of $\Kbar$
such that $\lim_{n\to \infty} h_{\F} (S_n) = 0$.
Then for all $v \in M_K$, the sequence $\{ S_n \}$
is pseudo-equidistributed with respect to the pair $(F_v,\lambda_v)$.
\end{theorem}

The proof is very similar to the proof of Theorem~\ref{thm2}, and we omit it.

\medskip

As an application of Theorem~\ref{thm2variant}, we consider the
``adelic Mandelbrot set''.  

\medskip

The classical Mandelbrot set $M = M_\infty$ can be defined as the locus of
points $c \in \C$ such that $0$ stays bounded under iteration of the
quadratic polynomial $z^2 + c$ (see \cite[Ch.~VIII]{Carleson-Gamelin} for more
details).
If we replace $\C$ by $\C_p$ for a prime number $p$ and take this as
the definition of the $p$-adic Mandelbrot set $M_p$, then it is easy
to see that  $M_p$ is just the closed unit disc in $\C_p$.

\medskip

It is well-known that the transfinite diameter of $M$ is 1.  We can
see this as follows.  Let $\phi_c(z) = z^2 + c$, and for $n\geq 0$ define 
$F_n := \{ c \in \C \; : \; \phi_c^n(0) \leq 2 \}$.  Then it is a
straightforward exercise to show 
(see \cite[Theorem VIII.1.2]{Carleson-Gamelin}) that 
$F_1 \supseteq F_2 \supseteq \cdots$ and 
$M = \cap_{n=0}^\infty F_n$.  
If we define $\psi_n(z) := \phi_z^n(0)$, then $\psi_n$ is a monic polynomial of degree
$2^{n-1}$ in $z$, and 
$F_n := \{ z \in \C \; : \; \psi_n(z) \leq 2 \}$.  

Therefore 
\[
c(M) = \lim_{n \to\infty} 2^{\frac{1}{2^{n-1}}} = 1.
\]

In particular, since we clearly have $c(M_p)=1$ for all primes $p$, we
have the product formula $c(\M) = 1$, where $\M = \{ M_v \}_{v \in
M_\Q}$ is the {\it adelic Mandelbrot set}.

\medskip

In view of the above comments, it is natural to attach to $\M$ an
adelic generalized Green's function $\g$ given by $\lambda_p(z) = \log^+|z|_p$
at the finite places of $\Q$, and given at the archimedean place by
\[
\lambda_\infty(z) := \lim_{n \to \infty} \frac{1}{2^{n-1}} \log^+|\frac{\psi_n(z)}{2}|.
\]

We then have a corresponding height function $h_{\M} :
\Qbar \to \R$.  It follows from the definitions that $h_\M(\alpha)=0$
if and only if the set $\{ \phi_\alpha^n(0) \}$ of orbits of $0$ under
$\phi_\alpha$ is finite.  
% \helpme{Is there a name for such points?}
% \fixme{From this point of view, $0$ is not special any more. 
% Fix a number $a \in K$ and ask the set of points $\alpha \in \Kbar$
% such that the set $\{\phi^n_{\alpha} ( a )\}$ of orbits of $a$ is 
% finite (i.e. $a$ is preperiodic under the iterates of $\phi_{\alpha}$. 
% Does this question also determine a \lq\lq canonical height\rq \rq on the 
% parameter space of quadratic maps?}

% We remark that the equilibrium measure $\mu_M$ for $M$ has support
% equal to $\partial M$.  This follows easily from Montel's theorem (see
% Theorem VIII.1.5 of \cite{Carleson-Gamelin}).

Applying Theorem~\ref{thm2variant} and Proposition~\ref{balancedprop}, we find:

\begin{theorem}
\label{thm2Mandelbrot}
Let $ \{ S_n \}$ be a sequence of distinct
finite Galois-stable subsets of $\Qbar$
such that $\lim_{n\to \infty} h_\M (S_n) = 0$.
Then for all $v \in M_\Q$, the sequence $\{ S_n \}$
is pseudo-equidistributed with respect to the pair $(M_v,\lambda_v)$.  In
particular, $\{ S_n \}$ is equidistributed with respect to the
equilibrium measure $\mu_M$ for $M$.
\end{theorem}

\section{Open questions}

We close the paper with a list of questions for further investigation.

\medskip

\begin{list1}

\item If $\phi \in \C_p[z]$ is a polynomial of degree $d \geq 2$, must
  the filled Julia set (resp. the Julia set) of $\phi$ be algebraically 
capacitable? More generally, under what condition the set $F$ determined
by a chain of lemniscates as defined above is algebraically capacitable?
\item Here is a related question (see \cite{Bezivin} for more
  details).  Let $\phi$ be as in the previous question, and suppose 
  that all periodic points of $\phi$ are repelling.  Is the filled
  Julia set (resp. the Julia set) of $\phi$ necessarily compact?  
%\item Is there a nonarchimedean analogue of Lemma~\ref{balancedlemma} if the filled Julia set of $\phi$ is compact?
%\item Is there a measure-theoretic equidistribution statement over
 % $\C_p$ analogous to Proposition~\ref{balancedprop}?
\item Can Theorem~\ref{thm2} be generalized to canonical heights with respect to
  more general dynamical systems on varieties?

% \item L.~DeMarco has recently introduced the notion of {\it
%     homogeneous capacity} for compact subsets of $\C^2$ 
% (see \cite{DeMarcothesis}).
%     An equivalent formulation, which makes sense over any valued
%     field, is as follows.  If $A \subseteq L^2$ is a bounded set,
%     define the {\it homogeneous transfinite diameter} of $A$ to be 
% \[
% c(A) := \lim_{n \to\infty} \sup_{z_1,\ldots,z_n \in A} \prod_{i \neq j} | z_i
% \wedge z_j |^
% {\frac{1}{n(n-1)}},
% \]
% where for $z=(x,y), z' = (x',y') \in L^2$ we put $z \wedge z' := xy' -
% x'y$.  (One can prove that this limit always exists and agrees with
% DeMarco's definition of homogeneous capacity for compact subsets of $\C^2$.)
% 
% If $\phi_1,\phi_2 \in L[z,w]$ are homogeneous polynomials of degree
% $d \geq 2$, we can define the filled Julia set $F$ of $\phi := (\phi_1,\phi_2) :
% L^2 \to L^2$ to be the set of points in $L^2$ which stay bounded under
% iteration of $\phi$.  Assume that $\phi$ is non-degenerate, in the
% sense that the resultant ${\rm Res}(\phi_1,\phi_2)$ is
% nonzero.
% Then DeMarco has proved in the case $L=\C$ that 
% \[
% c(F) = |{\rm Res}(\phi_1,\phi_2)|^{-\frac{1}{d(d-1)}}.
% \]
% Does the same formula holds for any algebraically closed local field $L$?
\end{list1}

\end{document}